\documentclass[11pt]{amsart}
\usepackage{amsfonts}
\usepackage{amsthm}
\usepackage{thmtools}

\usepackage[T1]{fontenc}

\usepackage{marvosym}
\usepackage{amsmath}
\usepackage{amssymb}
\usepackage{amsthm}
\usepackage{color}
\usepackage{graphicx}
\usepackage[applemac]{inputenc}
\usepackage{phonetic}
\usepackage{pb-diagram}
\usepackage[T1]{fontenc}
\usepackage{ stmaryrd }
\usepackage{url}

\font\msbmten msbm10
\def \Bbb#1{\hbox{\msbmten #1}}
\def\Uu{\Bbb U}
  \def\Nn{\Bbb N} \def\Rr {\Bbb R} 
\def\Ff{\Bbb F} 

 \def\Z{\mathbb Z}

\def\N{\mathbb N}

\newcommand\C{{\mathcal C}}

\newcommand\R{{\mathbb R}}

\def\Ind#1#2{#1\setbox0=\hbox{$#1x$}\kern\wd0\hbox to 0pt{\hss$#1\mid$\hss}     
\lower.9\ht0\hbox to 0pt{\hss$#1\smile$\hss}\kern\wd0}
\def\Notind#1#2{#1\setbox0=\hbox{$#1x$}\kern\wd0\hbox to 0pt{\mathchardef
\nn=12854\hss$#1\nn$\kern1.4\wd0\hss}\hbox to
0pt{\hss$#1\mid$\hss}\lower.9\ht0 \hbox to
0pt{\hss$#1\smile$\hss}\kern\wd0}

\def\ind{\mathop{\mathpalette\Ind{}}}                      
\def\nind{\mathop{\mathpalette\Notind{}}}

\def\dind{\mathop{\mathpalette\Ind{}}^{\delta}}  
\def\ndind{\mathop{\mathpalette\Notind{}}^{\delta} }

\def\la{\langle}
\def\ra{\rangle}
 
 \newcommand{\pte}[1]{\textcolor{red}{\textbf{#1}}}

\def\WOD{\mathop{\rm WOD}\nolimits}

\def\eq{\mathop{\rm eq}\nolimits}
\def\tp{\mathop{\rm tp}\nolimits}

\def\Fix{\mathop{\rm Fix}\nolimits}

\def\Aut{\mathop{\rm Aut}\nolimits}

\def\Th{\mathop{\rm Th}\nolimits}

\def\inf{\mathop{\rm inf}\nolimits}
\def\st{\mathop{\rm st}\nolimits}
\def\log{\mathop{\rm log}\nolimits}
\def\MD{\mathop{\rm MD}\nolimits}
\def\Max{\mathop{\rm max}\nolimits}
\def\Min{\mathop{\rm min}\nolimits}
\def\acl{\mathop{\rm acl}\nolimits}
\def\DC{\mathop{\rm DC}\nolimits}
\def\cl{\mathop{\rm cl}\nolimits}
\def\FMV{\mathop{\rm FMV}\nolimits}
\def\WOD{\mathop{\rm WOD}\nolimits}

\newtheorem{theorem}{Theorem}[subsection]

\newtheorem{lemma}[theorem]{Lemma}
\newtheorem{proposition}[theorem]{Proposition}
\newtheorem{conjecture}[theorem]{Conjecture}
\newtheorem{corollary}[theorem]{Corollary}
\newtheorem{problem}[theorem]{Problem}
\newtheorem{question}[theorem]{Question}

\newtheorem{fact}[theorem]{Fact}

\theoremstyle{definition}
\newtheorem{defi}[theorem]{Definition}
\newtheorem{remark}[theorem]{Remark}
\newtheorem{obs}[theorem]{Observation}

\newtheorem{example}[theorem]{Example}

\newcommand{\bp}{\begin{proposition}}
\newcommand{\ep}{\end{proposition}}
\newcommand{\bd}{\begin{defi}}
\newcommand{\ed}{\end{defi}}
\newcommand{\bej}{\begin{example}}
\newcommand{\eej}{\end{example}}
\newcommand{\bl}{\begin{lemma}}
\newcommand{\el}{\end{lemma}}
\newcommand{\bh}{\begin{fact}}
\newcommand{\eh}{\end{fact}}
\newcommand{\bpreg}{\begin{question}}
\newcommand{\epreg}{\end{question}}
\newcommand{\bo}{\begin{obs}}
\newcommand{\eo}{\end{obs}}
\newcommand{\bcon}{\begin{conjecture}}
\newcommand{\econ}{\end{conjecture}}
\newcommand{\brmk}{\begin{remark}}
\newcommand{\ermk}{\end{remark}}
\newcommand{\bc}{\begin{corollary}}
\newcommand{\ec}{\end{corollary}}
\newcommand{\be}{\begin{example}}
\newcommand{\ee}{\end{example}}
\newcommand{\bt}{\begin{theorem}}
\newcommand{\et}{\end{theorem}}

\newcommand{\bdem}{\begin{proof}}
\newcommand{\edem}{\end{proof}}
\newcommand{\benum}{\begin{enumerate}}
\newcommand{\eenum}{\end{enumerate}}
\newcommand{\bitem}{\begin{itemize}}
\newcommand{\eitem}{\end{itemize}}
\newcommand{\ov}[1]{\overline{#1}}

\begin{document}

\title{Pseudofinite structures and simplicity}
 
 \author{Dar\'io Garc\'ia}
\address{Dar\'io Garc\'ia  \\
Departamento de Matem\'aticas \\
Universidad de los Andes \\
Cra 1 No. 18A-10, Edificio H \\
Bogot\'a, 111711 \\
Colombia}

\email{dagarciar@gmail.com}

\author{Dugald Macpherson}
\address{Dugald Macpherson \\
School of Mathematics \\
University of Leeds  \\
Leeds LS2 9JT \\
UK}

\email{H.D.MacPherson@leeds.ac.uk}

\author{Charles Steinhorn}
\address{ Charles Steinhorn\\ Vassar College \\Poughkeepsie, NY
USA}
\email{steinhorn@vassar.edu}

\thanks{This paper is in part based on work supported by the U.S. National Science Foundation under Grant No. 0932078 000 while the authors were in residence at the Mathematical Sciences Research Institute in Berkeley, California during the Spring 2014 semester. Macpherson and Steinhorn were supported by EPSRC grant EP/K020692/1.}

\maketitle

{\bf Abstract.} We explore a notion of {\em pseudofinite dimension}, introduced by Hrushovski and Wagner, on an infinite ultraproduct of finite structures. Certain conditions on pseudofinite dimension are identified that guarantee simplicity or supersimplicity of the underlying theory, and  that a drop in pseudofinite dimension is equivalent  to forking. Under a suitable assumption, a measure-theoretic condition is shown to be equivalent to local stability. Many examples are explored, including vector spaces over finite fields viewed as 2-sorted finite structures, and homocyclic groups. Connections are made to products of sets in finite groups, in particular to word maps, and a generalization of Tao's algebraic regularity lemma is noted.

\section{Introduction}\label{qfd}

We investigate a notion of pseudofinite dimension (called {\em quasifinite dimension} in \cite{h1}) that was introduced in \cite{hw},  and applied by Hrushovski to approximate subgroups in \cite{h1} with tantalizing further directions suggested in \cite{Hru}. For an ultraproduct of finite structures,  the {\em pseudofinite dimension} $\delta(X)$ of a definable set $X$ is defined. It takes values in a quotient of the non-standard reals, rather than  in the positive integers or the ordinals, as holds for more standard model-theoretic dimensions and ranks. By taking an infimum in an appropriate completion, pseudofinite dimension is also defined for types. Given a definable set $X$, one obtains also a measure on the collection of its definable subsets which takes value 0 on strictly lower-dimensional subsets.

In the above papers, the main emphasis is on (pseudo)finite substructures of a given infinite structure such as an algebraically closed field or a simple algebraic group. A highlight is an abstract model-theoretic version of the Larsen-Pink inequality from \cite{lp}, linking e.g. pseudofinite dimension to Zariski dimension. 

Here, we explore general conditions on pseudofinite dimension that ensure simplicity, or stability, of the underlying theory, and yield a clear link between pseudofinite dimension and model-theoretic forking. Key conditions, introduced formally in Definition~\ref{extrahyp}, are (A), (SA), and $(\DC_L)$. These are all conditions on a pseudofinite structure $M$ that is an ultraproduct of finite structures. Roughly, (A) states that for any formula $\phi(\bar{x},\bar{y})$
the pseudofinite dimension of a consistent set of positive $\phi$-instances (a partial positive $\phi$-type) is obtained by a finite conjunction; it has a strengthening $A^*$, where positivity is not required.  The global version of (A), for an arbitrary partial type, is (SA). The condition $(\DC_L)$ roughly asserts that given an $L\/$-formula $\phi(\bar{x},\bar{y})$, the relation
$\delta(\phi(M^r,\bar{y}))<\delta(\phi(M^r,\bar{y}'))$ is definable by an $L\/$-formula $\psi(\bar{y},\bar{y}')$. We write $\bar{a}\dind_C B$ if $\delta(\tp(\bar{a}/B\cup C))=\delta(\tp(\bar{a}/C))$.

Our main results are as follows, all proved in Section~\ref{simplicityforking}. The assumptions are on an ultraproduct $M$ of a class of finite structures, though the conditions (A), (SA), and $(\DC_L)$ also make reference to a second sort in which counting takes place.


\begin{restatable}{thrm}{lowsimple}
\label{thm:lowsimple}
Assume that (A) holds. Then $\Th(M)$ is simple and low.
\end{restatable}



\begin{restatable}{thrm}{supersimple}
\label{thm:supersimple}
Assume (SA).
Then $\Th(M)$ is supersimple.
\end{restatable}


\begin{restatable}{thrm}{forking}
\label{thm:forking}
 Assume (SA) and $(\DC_L)$. Then for any $\bar{a},A,B$ in $M$ with $A,B$ countable,
$$\bar{a}\ind_A B\Leftrightarrow \bar{a}\dind_A B.$$
\end{restatable}

Working under a local version of $A^*$, we also characterize in Proposition~\ref{stablecrit} when a formula is stable, in terms of 
$\delta$ and a measure $\mu_D$ defined in Section~\ref{firstpreliminaries}.

In addition to examples constructed specifically  to delimit the conditions, some algebraically natural classes of examples are discussed in Section~\ref{examples}. 
These include ultraproducts of asymptotic classes of finite structures (from \cite{ms} and \cite{elwes}); the latter include the family of finite fields, every family of finite simple groups of fixed Lie type, and, for any smoothly approximable structure, an appropriately chosen family of envelopes. Ultraproducts of asymptotic classes satisfy $(\DC_L)$ and (SA). These conditions also hold for a (2-sorted) ultraproduct of finite vector spaces over finite fields, where in the ultraproduct the vector space dimension and the field are both infinite; unlike with asymptotic classes, the SU rank here is infinite in the vector space sort. An analogue in this setting (Theorem~\ref{cdmvs} below) is given of the main theorem of \cite{cdm}; we view this as a first example of a multi-sorted and infinite-rank enrichment of the notion of asymptotic class. We also give a uniformity result for exact (rather than asymptotic) cardinalities of definable sets in finite homocyclic groups (direct sums of isomorphic cyclic $p$-groups); see Proposition~\ref{abelianexact}.
This yields
an example  of an abelian group satisfying (A) and $(\DC_L)$ but not (SA). The second and third author, in conjunction with W. Anscombe and D. Wood, have work in progress on a rather flexible multi-sorted generalization of asymptotic classes, which should yield many more examples.

The paper is organized as follows. The framework, basic definitions, and easy observations around them are given in Section~\ref{preliminary}. In Section~\ref{deltavsforking} the main theorems stated above are proved. Section~\ref{examples} focuses on  examples, including the proof of Theorem~\ref{cdmvs}. 
Section~\ref{furtherprop} consists of some potentially useful technical results involving our conditions. These include the transfer of the conditions to $M^{\eq}$, one-variable criteria for the conditions, and results concerning a pregeometry defined via pseudofinite dimension (Propositions~\ref{geom} and~\ref{geom2}). Finally, in Section~\ref{applications}, we consider two possible lines of application of the results and framework. We note (Theorem~\ref{tao}) that the Pillay-Starchenko generalization of Tao's Algebraic Regularity Lemma holds under the assumptions (SA), $(\DC_L)$, and a corresponding definability condition for measure, and we also consider a possible application  of our results to pseudofinite groups (Theorem~\ref{wordmaps}). Open problems are mentioned throughout the paper. 

{\em Acknowledgement.} We thank William Anscombe for several very helpful conversations.

\section{Preliminaries}\label{preliminary}
\subsection{Pseudofinite structures, dimension, and measures}\label{firstpreliminaries}

We adopt the context of \cite[Section 5]{h1}, which extends \cite{hw}. We summarize it briefly.

We fix a countable first order language $L$, and consider $L$-formulas $\phi(\bar{x},\bar{y})$,  with the convention that $l(\bar{x})=r$ and $l(\bar{y})=s$. 
An {\em instance} of $\phi$ 
(in an $L$-structure $M$) is a formula $\phi(\bar{x},\bar{a})$, for $\bar{a}\in M^s$.
 A {\em $\phi$-formula} is a Boolean combination of instances of $\phi$. Parameter sets, always countable in this paper, are usually denoted by $A,B$, possibly with subscripts. Given a parameter set 
$B$ contained in a model $M$ of a theory $T$, we write $S_r(B)$ for the space of $r$-types of $T$ over $B$. A {\em partial positive $\phi$-type} is a  set of formulas
$\{\phi(x,\bar{a}):\bar{a}\in M^s\}$ that is consistent with $\Th(M,m)_{m\in M}$. A {\em $\phi$-type} over a parameter set is a consistent set of $\phi$-formulas over that set.

Let $\C$ be a class of finite $L$-structures. It is possible to extend $L$ to a language 
$L^+$ that includes a sort {\bf D} carrying the $L$-structure, a sort {\bf OF} carrying the language of ordered fields, and maps taking finite definable sets 
to their cardinalities. Formally, for each $L$-formula $\phi(\bar{x},\bar{y})$ with $\bar{x}=(x_1,\ldots,x_r)$ and $\bar{y}=(y_1,\ldots,y_s)$, the language  $L^+$ has a function $f_\phi:\textbf{D}^s \rightarrow \textbf{OF}$. 
Each structure $M_i\in \C$ gives a 2-sorted structure $K_i$ in $L^+$, the second sort 
being a copy of $(\Rr,<,+,\cdot,-,0,1,\log)$, and for $\phi(\bar{x},\bar{y})$ as above  and $\bar{a}\in M_i^s$, we put $f_\phi(\bar{a})=|\phi(M_i^r,\bar{a})|$. Let $\C^+:=\{K_i:M_i\in \C\}$ and denote by $K^*$ an ultraproduct of the members of $\C^+$ in $L^+$, over a non-principal ultrafilter $\mathcal{U}$. Here $K^{*}$ is 2-sorted with a sort {\bf OF} consisting
 of a non-archimedean real closed field $\R^*$, and a sort {\bf D} consisting of the ultraproduct of 
the $M_i$. Let $M$ denote the induced structure of $K^*$ in the {\bf D}-sort in the language $L$, and put $T:=\Th(M)$ and $T_M:=\Th(M, (m)_{m\in M})$. We usually refer just to the structure $M$ and theory $T$, without explicit reference to the 2-sorted context.

This construction might be seen as an explicit construction of the Cardinality Comparison Quantifiers $CCx$ defined by Hrushovski in \cite{Hru}. Indeed, with the notation used in  \cite{Hru}, we obtain that
\[M\models (CCx\,\phi)(b,b') \Leftrightarrow \text{$|\phi(M_i,b_i)|\geq|\phi(M_i,b'_i)|$ for $\mathcal{U}$-almost all $i$}\Leftrightarrow K^*\models f_\phi(b)\geq f_\phi(b').\]

Let $C$ be the convex hull of the integers $\mathbb Z$ in $\R^*$. Note that $C$ is a convex subgroup of $\R^*$. If $X$ is an $L$-definable set in $M$, say $X=\phi(M^r,\bar{a})$, we write $|X|$ for $f_\phi(\bar{a})$; observe that this is well-defined. Then
 define the \emph{pseudofinite dimension of $X$} to be $\delta(X):=\log|X|+C \in \R^*/C$. For definable sets $X,Y$ in $M$ we have 
 $$\delta(X)=\delta(Y) \mbox{~if and only if\ } 
\frac{1}{n}\leq \frac{|X|}{|Y|}\leq n \mbox{\ for some $n \in \Nn^{>0}$.}$$ 
Indeed, if $|X|\geq |Y|$, then
$$\delta(X)=\delta(Y)\Leftrightarrow \log|X|-\log|Y|\in C\Leftrightarrow \log(|X|/|Y|)\in C\Leftrightarrow |X|/|Y|\in C.$$ 
We write $\delta(\phi(\bar{x},\bar{b}))$ for $\delta(X)$ where $X$ is the set
 defined by the formula (with parameters) $\phi(\bar{x},\bar{b})$. In general the map $\bar{a}\mapsto \delta(\phi(\bar{x},\bar{a}))$ is not definable even in $L^+$, since $C$ and hence $\R^*/C$ are not definable.
 
 The map $\delta$ is extended to infinitely definable sets in \cite{h1}. For $\epsilon\in \R^*$, chosen sufficiently large and 
 with $\epsilon>C$, put 
$$V_0=V_0(\epsilon):= \{a\in \R^*/C: -n\epsilon+C \leq a\leq n\epsilon+C,n\in \Nn\}.$$ 
Let $V=V(\epsilon)$ be the set  of cuts in $V_0$, i.e., nonempty subsets bounded above and closed 
downwards. Then $V$ is a semigroup under set addition, and $V_0$ is identified with its image in $V$.
For a $\bigwedge$-definable set $X$, define
$$\delta(X):=\inf\{\delta(D): D\supset X, D \mbox{~definable}\},$$
the infimum evaluated in $V(\epsilon)$ for sufficiently large $\epsilon$.  
Given $B\subset M$ and a tuple $\bar{a}$ from $M$,
$\delta(\bar{a}/B)$ denotes $\delta(\tp(\bar{a}/B))$, and
$\delta^\phi(\bar{a}/B)$ denotes $\delta(\tp^\phi(\bar{a}/B))$, that is,  
the dimension of the corresponding $\phi$-type.

Hrushovski describes in  \cite{h1,Hru} different pseudofinite dimensions that can be obtained using the same construction of logarithms and taking quotients by different convex subgroups of $\R^*$. We focus in this paper on the pseudofinite dimension described above, also referred to as $\delta=\delta_{\text{fin}}$ in \cite{Hru}.
 We call $\delta$ a \emph{pseudofinite dimension} because its definition takes place in a pseudofinite context---an ultraproduct of finite structures---and it has some of the properties that we would expect for a dimension operator. For instance, the following are noted in \cite{hw} and \cite{h1}.

\bl \label{basics} \benum
\item[(i)]  $\delta(\emptyset)=-\infty$, and $\delta(X)=0$ for any finite definable set $X$.

\item[(ii)] If $X_1,X_2$ are $\bigwedge$-definable, then $\delta(X_1 \cup X_2) =\max\{\delta(X_1), \delta(X_2)\}$.

\item[(iii)] If $X_1,X_2$ are $\bigwedge$-definable, then $\delta(X_1 \times X_2) =\delta(X_1) +\delta(X_2)$.

\item[(iv)] If $(\alpha_n), (\beta_n)$ are descending sequences of cuts in $V_0$,
then $\inf_n(\alpha_n+\beta_n)=\inf_n\alpha_n+\inf_n \beta_n$.

\item[(v)] If $\alpha,\alpha',\beta,\beta'\in V$ with $\alpha<\alpha'$ and $\beta<\beta'$ 
then $\alpha+\beta<\alpha'+\beta'$.

\item[(vi)] If $X=\bigcap X_n$ with $X_1 \supset X_2 \supset\ldots$ all $\bigwedge$-definable, then
$\delta(X)=\inf_n\delta(X_n)$. 

\item[(vii)] If $X$ is $\bigwedge$-definable, $f$ is a definable map, $\gamma \in V_0$, and 
$\delta(f^{-1}(a)\cap X)\leq \gamma$ for all $a$, then $\delta(X)\leq \delta(f(X))+\gamma$.

\eenum
\el

For a definable subset $D$ of $M$, there is a finitely additive real-valued probability measure $\mu_D$ on the set 
of definable subsets $X$ of $D$ given by $\mu_D(X):=\st\left(\lim_{i\to \mathcal{U}} \frac{|X(M_i)\cap D(M_i)|}{|D(M_i)|}\right)$, where $\st(\cdot)$ is the standard part map. This measure can be extended  to  a countably-additive probability measure on the 
$\sigma$-algebra generated by the definable subsets of $D$, and thus is defined on the $\bigwedge$-definable subsets of $D$. 
This measure combines with $\delta$ in our characterization of when a formula is stable, Proposition~\ref{stablecrit}.

\subsection{Conditions on the pseudofinite dimension}\label{dimensionconditions}

 We investigate the following hypotheses on $\delta$. Throughout, we work in the context described in Section~\ref{firstpreliminaries}, with $M$ the $L$-structure induced on the sort ${\bf D}$ by an infinite ultraproduct of  finite $L^+$-structures.  As usual, let $\phi(\bar{x},\bar{y})$ be an $L$-formula, with $l(\bar{x})=r$ and $l(\bar{y})=s$. For a type $p$ over $B$ and $B_0\subset B$, we let $p|B_0$ denote the restriction of $p$ to $B_0$.

\bd \label{extrahyp}
\benum
\item \emph{Attainability} ($A_\phi$). There is no sequence $(p_i:i \in \omega)$ of finite partial positive $\phi$-types such that  $p_i \subseteq p_{i+1}$ (as sets of formulas) and $\delta(p_i)>\delta(p_{i+1})$ for each $i\in \omega$. We denote by ($A_\phi^*$) the corresponding (stronger) condition, where the above is assumed for all increasing sequences of finite partial $\phi$-types, not necessarily positive.

\item \emph{Strong Attainability} ($SA$). For each partial type $p(\bar{x})$ over a parameter set $B$, there is a finite subtype $p_0$ of $p$
 such that $\delta(p(\bar{x}))=\delta(p_0(\bar{x}))$. 

\item ($SA^-$). If $(B_i)_{i\in\omega}$ is a sequence of countable parameter sets with $B_i\subset B_{i+1}$ for each $i$, $B:=\bigcup_{i\in\omega}B_i$ and $p\in S(B)$, then there is $j\in\omega$ such that $\delta(p|B_j)=\delta(p|B_i)$ for all $i\geq j$. 

\item \emph{Weak Order Definability} ($\WOD_\phi$). There is $n=n_\phi\in \Nn$ such 
that for all $\bar{a},\bar{b}\in M^s$,
$$\delta(\phi(\bar{x},\bar{a}))=\delta(\phi(\bar{x},\bar{b}))\Leftrightarrow \frac{1}{n}<\frac{|\phi(\bar{x},\bar{a})|}{|\phi(\bar{x},\bar{b})|}\leq n.$$

\item \emph{Dimension Comparison in $L^+$} ($\DC_{L^+}$).  For all formulas $\phi(\bar{x},\bar{y})$ and $\psi(\bar{x},\bar{z})$ (with $t=l(\bar{z})$), there is an $L^+$-formula
$\chi_{\phi,\psi}(\bar{y},\bar{z})$ such that for all $\bar{a}\in M^s$ and $\bar{b}\in M^t$,
$$\chi_{\phi,\psi}(\bar{a},\bar{b})\Leftrightarrow \delta_{\bar{x}}(\phi(\bar{x},\bar{a}))\leq \delta_{\bar{x}}(\psi(\bar{x},\bar{b})).$$

\item \emph{Dimension Comparison in $L$} ($\DC_L$). This is as for ($\DC_{L^+}$), except that the formula  $\chi_{\phi,\psi}$ can be chosen in $L$.

\item \emph{Finitely Many Values} ($\FMV_\phi$). There is a finite set $\{\delta_1,\ldots,\delta_k\}$ such that for each $\bar{b}\in M^s$ there is $i\in\{1,\ldots,k\}$ with
$\delta(\phi(M^r,\bar{b}))=\delta_i$.
\eenum
\ed

The conditions ($A_\phi$), ($A_\phi^*$), ($\WOD_\phi$) and ($\FMV_\phi$) have global  versions ($A$), ($A^*$), ($\WOD$) and 
($\FMV$), where they are assumed to hold for all $\phi$ (with $k$ and the $\delta_i$ in $(\FMV)$ dependent on $\phi$). They can also be formulated for finite sets $\Delta$ of formulas, and it would be helpful to know to what extent they---in particular, local versions of ($A$)---are preserved under taking Boolean combinations of formulas.

We conclude this section with some easy observations about these conditions. Note that, trivially, ($SA$) implies ($SA^-$), and observe that 
 $$\mbox{($SA$)\ \ }\Rightarrow \mbox{\ \ ($A_\phi^*$) for all $\phi$\ \ } \Rightarrow\mbox{\ \ ($A$).}$$
We also have 
  
 \begin{lemma} \label{A and A*}
(i) For every formula $\phi(\bar{x},\bar{y})$, the conditions 
$((A_\phi)\wedge (A_{\neg \phi}) \wedge (A_{\phi(\bar{x},\bar{y}_1)\wedge \neg \phi(\bar{x},\bar{y}_2)}))$ and ($A_\phi^*$) are equivalent.

(ii) The conditions ($A$) and ($A^*$) are equivalent.
\end{lemma}

\bdem
(i) The implication ($\Leftarrow$) is immediate. For the direction ($\Rightarrow$), suppose ($A_\phi^*$) fails, witnessed by an infinite decreasing sequence of finite partial $\phi$-types $p_1\subset p_2\subset\ldots$ with $\delta(p_1)>\delta(p_2)>\ldots$. 
If the $p_i$ do not involve $\neg\phi$ then ($A_\phi$) fails, and if $\phi$ does not occur positively in the $p_i$ then 
($A_{\neg \phi}$) fails. Otherwise, we may suppose both $\phi$ and $\neg\phi$ occur in $p_1$. It is now easy to construct an increasing sequence of positive
 $(\phi(\bar{x},\bar{y}_1)\wedge \neg \phi(\bar{x},\bar{y}_2))$-types with strictly decreasing $\delta$-values, repeating some parameters if needed.

(ii) This is immediate from (i). 
\edem

We shall make little use of ($\WOD$) and ($\DC_{L^+}$), but we note:

\begin{lemma} \label{WOD}
The conditions ($\WOD$) and  ($\DC_{L^+}$) are equivalent.
\end{lemma}
\bdem
$\mbox{($\DC_{L^+}$)\ }\Rightarrow \mbox{\ ($\WOD $)}$.
By compactness and $\omega_1$-saturation of $M$, ($\DC_{L^+}$) implies ($\WOD$). Indeed, assume ($\DC_{L^+})$, and
 that for all $n\in \Nn$ there are $\bar{a}_n,\bar{b}_n \in M^s$ with 
$$\delta(\phi(\bar{x},\bar{a}_n))=\delta(\phi(\bar{x},\bar{b}_n))\text{\ \ and\ \ }
 \frac{|\phi(\bar{x},\bar{a}_n)|}{|\phi(\bar{x},\bar{b}_n)|}<\frac{1}{n}.$$
 Then the set of $L^+$-formulas
 $$\{\chi_{\phi,\phi}(\bar{y}_1,\bar{y}_2), \chi_{\phi,\phi}(\bar{y}_2,\bar{y}_1)\}\cup 
\left\{\frac{|\phi(\bar{x},\bar{y}_1)|}{|\phi(\bar{x},\bar{y_2})|}<\frac{1}{n}:n\in\omega\right\}$$
is consistent. This is impossible, by our initial remarks on pseudofinite dimension.

$\mbox{($\WOD$)\ }\Rightarrow \mbox{\ ($\DC_{L^+}$)}$.
Fix formulas $\phi(\bar{x},\bar{y})$ and $\psi(\bar{x},\bar{z})$. We may regard the union of the family of sets defined by $\phi$ and the family of sets defined by $\psi$, as a single uniformly definable family of sets, defined e.g. by the formula $\rho(\bar{x}, \bar{y}\bar{z}ww_1w_2)$ which has form
$$((w=w_1\wedge w_1\neq w_2)\to \phi(\bar{x},\bar{y})) \wedge (w=w_2 \wedge w_1\neq w_2) \to \psi(\bar{x},\bar{z}))$$
$$\wedge ((w_1=w_2 \vee (w\neq w_1 \wedge w\neq w_2)) \to (\neg \phi(\bar{x},\bar{y}) 
\wedge \neg \psi(\bar{x},\bar{z})).$$
By ($\WOD$), there is a number $n_\rho$ associated with $\rho$. Now
$\delta(\phi(\bar{x},\bar{a}))\leq \delta(\psi(\bar{x},\bar{b}))$ holds if and only if 
$$\exists w_1,w_2(w_1\neq w_2 \wedge |\rho(\bar{x},\bar{a}\bar{b}w_1w_1w_2)|\leq n_\rho|\rho(\bar{x},\bar{a}\bar{b}w_2w_1w_2)|).$$
\edem

We also note that under ($\WOD$) the following strengthening of Lemma~\ref{basics}(vii) holds for definable sets. We omit the proof, an easy counting argument in finite structures.

\begin{lemma}\label{WODapp}
(i) Let $f:X\to Y$ be a definable map in $M$ between definable sets $X,Y$, and suppose that there is a positive integer $n$ such that for all $\bar{a},\bar{b}\in Y$,
$\frac{1}{n}|f^{-1}(\bar{b})|\leq |f^{-1}(\bar{a})|\leq n|f^{-1}(\bar{b})|$ (non-standard cardinalities). Then for all $\bar{a}\in Y$, we have 
$\delta(X)=\delta(Y)+\delta(f^{-1}(\bar{a}))$.

(ii) Assume $M$ satisfies ($\WOD$). Let $X$ be a definable set in $M$, $f$ a definable map, and suppose that $\delta(f^{-1}(\bar{a}))=\gamma$ for all $\bar{a}\in f(X)$. 
Then $\delta(X)=\delta(f(X))+\gamma$.
\end{lemma}

\begin{lemma}\label{Acontent}
Assume that ($A_\phi$) holds. Then  there is $m=m_\phi\in \omega$ such that there do not exist $\bar{a}_1,\ldots,\bar{a}_m\in M^s$ so that if $p_i:=\{\phi(\bar{x},\bar{a_j}):j\leq i\}$ then $p_i$ is consistent and 
$\delta(p_1)>\delta(p_2)>\ldots>\delta(p_m)$.
\end{lemma}

\bdem Suppose not. Then for every $N\in \omega$ there are $\bar{a}_1,\ldots,\bar{a}_N$ such that 
$$|\phi(\bar{x},\bar{a}_1)\wedge \ldots \wedge \phi(\bar{x},\bar{a}_i)|>N|\phi(\bar{x},\bar{a}_1)\wedge \ldots \wedge \phi(\bar{x},\bar{a}_{i+1})|$$
for each $i=1,\ldots, N$. Each such statement is a partial type in $L^+$. It follows by compactness and $\omega_1$-saturation of $K^{*}$ that there are
$\bar{a}_i\in M^s$ for all $i>0$ such that for each $i,N \in \omega$, we have
$$|\phi(\bar{x},\bar{a}_1)\wedge \ldots \wedge \phi(\bar{x},\bar{a}_i)|>N|\phi(\bar{x},\bar{a}_1)\wedge \ldots \wedge \phi(\bar{x},\bar{a}_{i+1})|.$$
Putting $p_i:=\{\phi(\bar{x},\bar{a}_j): j\leq i\}$, we obtain an infinite  
sequence $(p_i)_{i>0}$ of finite partial positive $\phi$-types with $\delta(p_i)>\delta(p_{i+1})$ for each $i$, contrary to ($A$). \edem


\bl Assume ($SA$) holds. Then there is no sequence of definable sets $(S_n:n<\omega)$ such that $S_{n+1}\subseteq S_n$ and $\delta(S_{n+1})<\delta(S_n)$ for each $n<\omega$.
\el
\bdem Suppose $S_i$ is defined by the formula $\phi_i(\ov{x};\ov{b}_i)$ and consider the partial type \[p:=\{\phi_i(\ov{x};\ov{b}_i):i<\omega\}\] By (SA), there is $m<\omega$ such that $\delta(p)=\delta(\phi_m(\ov{x},\ov{b}_m))$, which is impossible because $\delta(p)\leq \delta(S_{m+1})<\delta(S_m)$. \edem

\begin{lemma} \label{SOP} 
Assume ($\DC_L$), and suppose ($\FMV_\phi$) fails for some $\phi$. Then $T$ has the strict order property, so in particular is not simple. 
\end{lemma}

\bdem
There is a preorder with an infinite chain definable on $M^s$, where we put
$\bar{y}\leq \bar{y}'\Leftrightarrow \delta(\phi(\bar{x},\bar{y}))\leq \delta(\phi(\bar{x},\bar{y}'))$. 
\edem

\begin{proposition} \label{sopei}
(i) Assume $\Th(M)$ has the strict order property. Then  ($\FMV$) fails.

(ii) If $M$ has ($\DC_L$) and $\Th(M)$ has elimination of imaginaries (EI), then
($\FMV$) holds, and so
$\Th(M)$ does not have the strict order property.
\end{proposition}

\bdem (i) Let $\psi(\bar{u},\bar{v})$ be a formula defining a preorder $\preceq$ on $M^t$ with an infinite chain.  Find by $\omega_1$-saturation an infinite $L^+$-indiscernible sequence $(\bar{a}_i:i\in\omega)$ in $M^t$, with $\bar{a}_i\prec \bar{a}_j$ if and only if $i<j$, for $i,j\in \omega$. Let $\chi(\bar{x},\bar{u}\bar{v})$ express $\bar{u}\prec\bar{x}\prec\bar{v}$. 
Clearly, for any $n>0$, $f_\chi(\bar{a}_0\bar{a}_{n})\geq nf_\chi(\bar{a}_0\bar{a}_{1})$, where $f_\chi$ is the $L^+$-function symbol corresponding to $\chi$; indeed, 
$\chi(M^t,\bar{a}_0\bar{a}_n)\supseteq \bigcup_{i=0}^{n-1}\chi(M^t,\bar{a}_i\bar{a}_{i+1})$ and $f_\chi(\bar{a}_i\bar{a}_{i+1})=f_\chi(\bar{a}_0\bar{a}_1)$ for each $i$.  Hence by $L^+$-indiscernibility
$f_\chi(\bar{a}_0\bar{a}_2) \geq nf_\chi(\bar{a}_0\bar{a}_1)$ for each $n$, and so $\delta(\chi(M^t,\bar{a}_0\bar{a}_2))>\delta(\chi(M^t,\bar{a}_0\bar{a}_1))$. It follows by compactness that the set $\{\delta(\chi(M^t,\bar{a}_0\bar{a}_i)):i>0\}$ is infinite.

(ii) Assume ($\DC_L$) and EI, and suppose for a contradiction that $\{\delta(\phi(M^r,\bar{a})):\bar{a}\in M^s\}$ is infinite. Let $\psi(\bar{u},\bar{v})$ express that $\delta(\phi(M^r,\bar{u}))\leq\delta(\phi(M^r,\bar{v}))$. Then $\psi$ defines a pretotal order on $M^s$. Let $E$ be the equivalence relation defined by putting 
$E(\bar{u},\bar{v})\Leftrightarrow (\psi(\bar{u},\bar{v})\wedge \psi(\bar{v},\bar{u}))$. Then $\psi$ induces on $M^s/E$  an $\emptyset$-definable total order $<$. By elimination of imaginaries, for some $t$ there is an $\emptyset$-definable function $g:M^s\to M^t$ with $E(\bar{u},\bar{v})\Leftrightarrow  (g(\bar{u})=g(\bar{v}))$. There is then an $\emptyset$-definable total order $\prec$ on the infinite set $I:=g(M^s)$, given by $\bar{a}\prec \bar{b}$ if and only if $g^{-1}(\bar{a})<g^{-1}(\bar{b})$. Since $I$ is pseudofinite we may find a  sequence of subintervals $I\supset J_0\supset J_1\supset J_2\supset \ldots$ with $|J_i|=2|J_{i+1}|$ (non-standard cardinalities). Since intervals are uniformly definable, this contradicts ($\WOD$), and hence ($\DC_L$) by Lemma~\ref{WOD}. \edem

\noindent In Section~\ref{furtherprop} we say more about imaginaries, in particular about the preservation of our conditions when passing to $M^{\eq}$.

\begin{example} \label{earlyexample}\rm   It is easy to produce examples satisfying ($A$) or ($SA$) but without ($\FMV$). 
Let $L$ be a language with a single binary relation $E$, and let $M_k$ be an $L$-structure with $\displaystyle{\sum_{i=1}^k i^2}$ elements, in which $E$ is interpreted by an equivalence relation with a class of size $i^2$ for each $i=1,\ldots,k$. 
It can be shown that any non-principal ultraproduct $M$ satisfies ($SA$)  but not ($\DC_{L^+}$) or ($\FMV$). 
\end{example}

The following measure-theoretic result will be applied to the measures $\mu_D$ in the proof of Proposition~\ref{stablecrit} in Section~\ref{deltavsforking}. 

\begin{proposition} \label{kintersections} Let $X$ be a measure space with $\mu(X)=1$ and fix $0<\epsilon\leq \dfrac{1}{2}$. Let $\langle A_i:i<\omega\rangle$ be a sequence of measurable subsets of $X$ such that $\mu(A_i)\geq \epsilon$ for every $i$.
Then, for every $k<\omega$, there are $i_1<i_2<\ldots<i_k$ such that \[\mu\left(\bigcap_{j=1}^k A_{i_j}\right)\geq \epsilon^{3^{k-1}}\] 
\end{proposition}
\bdem The proof is by induction on $k$. 

$k=1$.\enspace By hypothesis we have $\mu(A_i)\geq \epsilon= \epsilon^{3^{1-1}}$. 

$k=2$.\enspace Assume the conclusion is false. Then $\mu(A_i\cap A_j)<\epsilon^{3^{2-1}}=\epsilon^3$ for all $i\neq j$. By the truncated Inclusion-Exclusion Principle we know that for every $N\in\N$,
\begin{align*}
1 &\geq \mu\left(\bigcup_{i=1}^N A_i\right)\geq \sum_{i=1}^N \mu(A_i) - \sum_{1\leq i<j\leq N} \mu(A_i\cap A_j)\geq N\epsilon - \dfrac{N(N-1)}{2}\epsilon^3. &(\dag)
\end{align*}
Define $f(x):=x\cdot \epsilon - \dfrac{x(x-1)}{2}\epsilon^3=-\dfrac{x^2}{2}\epsilon^3 + x\left(\epsilon+\dfrac{\epsilon^3}{2}\right).$

This function achieves its maximum value at $x_0=\dfrac{1}{\epsilon^2}+\dfrac{1}{2}$>0, and by taking any integer $N\in [x_0-1,x_0]$ we have that 
\begin{align*}
f(N)\geq f(x_0-1)&=\left(\dfrac{1}{\epsilon^2}-\dfrac{1}{2}\right)\epsilon - \dfrac{\left(\dfrac{1}{\epsilon^2}-\dfrac{1}{2}\right)\left(\dfrac{1}{\epsilon^2}-\dfrac{3}{2}\right)}{2} \cdot \epsilon^3\\
&=\dfrac{1}{2\epsilon}+\dfrac{\epsilon}{2}-\dfrac{3}{8}\epsilon^3\\
&\geq 1+\epsilon \left(\dfrac{1}{2}-\dfrac{3}{8}\epsilon^2\right) \hspace{1cm}\left[\text{because $\epsilon\leq \frac{1}{2}$}\right]\\
&>1.
\end{align*} contradicting (\dag). 

{\em Induction Step}.\enspace By the induction hypothesis, we can assume that there is a tuple $(i_1,\ldots,i_k)$ satisfying \[i_1<\ldots < i_k \text{\,\,\, and \,\,\,} \mu\left(\bigcap_{j=1}^k A_{i_j}\right)\geq \epsilon^{3^{k-1}}.\text{\hspace{1cm}} (*)\]

We claim that there are infinitely many such tuples. Indeed,
if not, take $\ell$ to be the maximum of all indices appearing in the tuples $(i_1,\ldots,i_k)$ that satisfy (*). Then, $\langle A_{j}:j\geq \ell+1\rangle$ is a sequence contradicting the induction hypothesis. 

Now let $\langle \alpha_j:j<\omega\rangle$ be an enumeration of all tuples satisfying (*) and put $\displaystyle{B_j=\bigcap_{i\in \alpha_j}A_i}$. By construction, $\mu(B_j)\geq \epsilon^{3^{k-1}}$ for all $j$.

By the $k=2$ case, there are indices $j_1\neq j_2$ such that \[\mu(B_{j_1}\cap B_{j_2})\geq (\epsilon^{3^{k-1}})^3=\epsilon^{3^{k-1}\cdot 3}=\epsilon^{3^k}\] where $j_1,j_2$ are indices corresponding to two different tuples $\alpha_{j_1}\neq \alpha_{j_2}$. In particular, there are (at least) $k+1$ indices $i_1<i_2<\cdots<i_k<i_{k+1}$ such that
\[\mu\left(\bigcap_{j=1}^{k+1}A_{i_j}\right)\geq \mu(B_{j_1}\cap B_{j_2})\geq \epsilon^{3^k}=\epsilon^{3^{(k+1)-1}} \hspace{1cm}\]  
\edem

\brmk\label{k=2} 1. In the $k=2$ case of the proof above, we actually find a number $N=N(\epsilon)$ such that if $A_1,\ldots,A_N$ have measure at least $\epsilon$, there are $1\leq i<j\leq N$ such that $\mu(A_i\cap A_j)\geq \epsilon^3$. 

2.   For the measures $\mu_D$, it is possible to provide a ``pseudofinite'' proof of Proposition \ref{kintersections}
that uses finite combinatorics and counting transferred to the ultraproduct via the functions $f_\phi$. 
\ermk

\section{Forking independence and $\delta$-independence} \label{deltavsforking}

\noindent Here we investigate $\delta$-independence, and prove Theorems~\ref{thm:lowsimple}--\ref{thm:forking} stated in the introduction. We also give a criterion for stability of a formula, Proposition~\ref{stablecrit}.

Our context is that of the whole paper: $\C$ is an infinite class of finite $L$-structures, $\C^+$ is the corresponding class of 
2-sorted $L^+$-structures in the sorts ${\bf D}$ and 
${\bf OF}$,  ${\mathcal U}$ is a non-principal ultrafilter on $\C^+$, and $K^*=(M^+,\R^*)$ is an ultraproduct of the structures in 
$\C^+$ with respect to ${\mathcal U}$. As usual, we denote by $M$ the reduct of $M^+$ to $L$, and by $T$ its theory. 
 
\subsection{Properties of $\delta$-independence.}

\bd \label{delta forking} Let $\bar{a}$ be  a tuple and $A, B$ be countable subsets of $M$. We say that $\bar{a}$ is \emph{$\delta$-independent of $B$ over $A$}, written $\displaystyle{a\dind_A B}$, if $\delta(\bar{a}/AB)=\delta(\bar{a}/A)$; here, as usual, $AB$ is an abbreviation for $A\cup B$. 
\ed

\begin{remark}\rm
With $\bar{a},B,A$ as in Definition~\ref{delta forking}, $\bar{a}\nind^\delta_A B$ if and only if there is a formula $\theta(\bar{x})\in\tp(\bar{a}/AB)$ such that for all $\psi(\bar{x})\in \tp(\bar{a}/A)$ we have $\delta(\theta(\bar{x}))<\delta(\psi(\bar{x}))$. The direction 
($\Rightarrow$) is immediate by Lemma~\ref{basics}(vi). The direction ($\Leftarrow$) requires a small compactness and saturation argument, using essentially that a cut in $\R^*/C$ cannot have a countably infinite initiality. More generally, if $p$ is a type over a countably infinite set and there is a formula $\theta$ with $\delta(p)=\delta(\theta)$, then there is $\psi\in p$ with $\delta(p)=\delta(\psi)$.
\end{remark}

We investigate below the extent to which  $\delta$-independence satisfies standard properties of non-forking. Throughout this section, we write $\ind$ for the usual forking independence relation.

\begin{lemma} \label{fiber} \emph{[Additivity]}
 Assume ($\DC_L$) and ($\FMV$), and let $A$ be a countable set of parameters from $M\models T$. Let $\bar{a}\in M^r,\bar{b}\in M^s$. Then $\delta(\bar{a}\bar{b}/A)=\delta(\bar{a}/A\bar{b})+\delta(\bar{b}/A)$. 
\end{lemma}

\bdem
Since $A$ is countable it plays no role in the proof and we thus may suppress parameters and suppose $A=\emptyset$. Let $(\phi_n(\bar{y}))_{n\in\omega}$ enumerate the formulas in $\tp(\bar{b})$ and $(\psi_n(\bar{x},\bar{b}))_{n\in\omega}$ enumerate those in $\tp(\bar{a}/\bar{b})$. We may suppose that for each $n\in \omega$, $\phi_{n+1}\to \phi_n$ and $\psi_{n+1} \to \psi_n$. Let $P$ be the set of realizations in $M$ of $\tp(\bar{a}\bar{b})$.
 
 Put $\epsilon_n :=\delta(\phi_n(\bar{y}))$ and $\gamma_n:=\delta(\psi_n(\bar{x},\bar{b}))$.
 For each $n$  there is $\rho_n(\bar{y}) $ over $\emptyset$ which expresses that $\delta(\{\bar{x}:\psi_n(\bar{x},\bar{y})\})=\gamma_n$. 
Indeed, by ($\FMV$), $\delta(\psi(M^r,\bar{b}'))$ takes finitely many values as $\bar{b}'$ varies, and by ($\DC_L$) we may compare these and thus express that the $j^{th}$ such value is taken.
 There is $f(n)\geq n$ such that
$\phi_{f(n)}(\bar{y}) \rightarrow \rho_n(\bar{y})$, and by refining the sequence $(\phi_n)_{n\in\omega}$ we may suppose that $\phi_n\to \rho_n$. 
Let $P_n$ be the set defined by $\phi_n(\bar{y}) \wedge \psi_n(\bar{x},\bar{y})$. 

We claim that  
 $\delta(P_n)=\epsilon_{n}+\gamma_n$ for each $n$. Indeed,
for a fixed $N>0$, by counting in finite structures we have
$$\frac{1}{N}|\psi_n(\bar{x},\bar{b})|\cdot |\phi_n(\bar{y})|\leq |P_n|\leq N|\psi_n(\bar{x},\bar{b})|.|\phi_n(\bar{y})|.$$
The claim now follows by taking logarithms and working modulo $C$ (see Lemma~\ref{WODapp}(i)).

Also, $P=\bigcap_n P_n$.
  By Lemma~\ref{basics}(vi), 
$$\delta(P)=\inf_n\delta(P_n)=\inf_n(\epsilon_n+\gamma_n).$$
Put $\epsilon:=\inf_n\epsilon_n =\delta(\tp(\bar{b}))$ and $\gamma:=\inf_n\gamma_n=\delta(\tp(\bar{a}/\bar{b}))$. Then by Lemma~\ref{basics}(iv),
$\epsilon+\gamma=\inf_n(\epsilon_n+\gamma_n)$, which equals $\delta(P)$, as required. 
\edem

\bp \label{d-properties} The following are properties of the $\delta$-independence relation:
\benum
\item[(i)] \emph{Existence}: Given countable sets $A\subset B$ and $p\in S_r(A)$ (for any $r\in \Nn$) there is
 $\bar{a}\models p$ with 
 $\bar{a}\dind_A B$.
\item[(ii)] \emph{Monotonicity and transitivity:}  If $A\subset D \subset B$, then 
$$\bar{a}\dind_A B\Leftrightarrow \left(\bar{a}\dind_A D \mbox{~and~} \bar{a}\dind_D B\right).$$
\item[(iii)] \emph{Finite character:} If $\displaystyle{a\nind^\delta_A B}$ then there is a finite subset $\ov{b}\subseteq B$ such that $\displaystyle {a\nind^\delta_{A}\ov{b}}$.
\eenum
\ep
 \bdem  \emph{Monotonicity and transitivity}.\enspace These follow immediately from the definitions. 
 
\noindent \emph{Existence}.\enspace This is easy to prove by compactness. We must show that given a partial type $q$ over $B$ and a
 formula $\phi(\bar{x},\bar{b})$ over $B$, if $\delta(q)=\delta_0$ then either 
 $\delta(q \cup \{\phi(\bar{x},\bar{b})\})=\delta_0$ or 
$\delta(q \cup \{\neg\phi(\bar{x},\bar{b})\})=\delta_0$. If this were to fail then there would be
 $\psi \in q$
such that $\delta(\psi\wedge\phi)<\delta_0$ and $\delta(\psi\wedge \neg \phi)<\delta_0$. Since
 $\delta(\psi)\geq \delta_0$, this contradicts Lemma~\ref{basics}(ii). 

\noindent\emph{Finite character}.\enspace Suppose that $\bar{a}\ndind_A B$. Then
 $\delta(\bar{a}/AB)<\delta(\bar{a}/A)$, so there is a  formula $\phi(\bar{x},\bar{b})$ over $B$ such that
 $\delta(\tp(\bar{a}/A) \cup \{\phi(\bar{x},\bar{b})\})<\delta(\bar{a}/A)$. Then $\bar{a}\ndind_A \bar{b}$.
\edem

\bp \label{d-properties+} Under the further assumptions listed below, we have:
\benum
\item[(iv)] \emph{Local character:} ($A$)\enspace For every $\bar{a}$ and $B\subseteq M$, there is a countable subset $A\subseteq B$ such that \[\displaystyle{\bar{a}\dind_A B}\]
\item[(v)] \emph{Invariance:} ($\DC_L$)\enspace If $\alpha\in \Aut(M)$, then $\displaystyle{\bar{a}\dind_A B \Leftrightarrow \alpha(\bar{a})\dind_{\alpha(A)}\alpha(B)}$
\item[(vi)] \emph{Symmetry:} ($\DC_L$)\enspace and $(\FMV)$) $\displaystyle{\bar{a}\dind_A \bar{b}}$ if and only if 
$\displaystyle{\bar{b}\dind_A \bar{a}}$ 
\item[(vii)] \emph{Algebraic closure:} ($\DC_L$)\enspace If $A\subseteq B\subseteq \acl^{\eq}(A)$, then $\delta(\bar{a}/A)=\delta(\bar{a}/B)$, where $\delta$ is defined in the natural way for formulas with parameters in $M^{\eq}$.
\eenum
\ep 
\bdem
\emph{Local character.}\enspace Let $p:=\tp(\bar{a}/B)$. By ($A$), for each $L$-formula $\phi(\bar{x},\bar{y})$ with $l(\bar{x})=l(\bar{a})$, there is a $\phi$-formula
$\psi_\phi(\bar{x},\bar{b_\phi})$ (a conjunction of $\phi$-instances) such that 
$\delta^\phi(\bar{a}/B)=\delta(\psi_\phi(\bar{x},\bar{b}_\phi))$. 
If $A$ is the collection of all elements in the tuples $\bar{b}_\phi$ as $\phi$ varies, then
$\bar{a}\dind_A B$ and $|A|\leq \aleph_0$.

\noindent \emph{Invariance}.\enspace Suppose $\bar{a}\dind_A B$ and $\alpha \in \Aut(M)$. For 
every $\phi(\bar{x},\bar{b})\in \tp(\bar{a}/AB)$ there is 
$\psi(\bar{x},\bar{c})\in \tp(\bar{a}/A)$ such that
$\delta(\psi(\bar{x},\bar{c}))\leq \delta(\phi(\bar{x},\bar{b}))$. By ($\DC_L$), we have $\chi_{\psi,\phi}(\bar{c},\bar{b})$, so
$\chi_{\psi,\phi}(\alpha(\bar{c}),\alpha(\bar{b}))$, noting that $\chi_{\psi,\phi}$ is an $L$-formula so is preserved by $\alpha$. Hence, 
$\delta(\psi(\bar{x},\alpha(\bar{c})))\leq \delta(\phi(\bar{x},\alpha(\bar{b})))$, and
 it follows that
$\alpha(\bar{a})\dind_{\alpha(A)} \alpha(B)$.

\noindent \emph{Symmetry}.\enspace  
It suffices to show that 
$\bar{a}\dind_A \bar{b}\Rightarrow \bar{b}\dind_A\bar{a}$,
so we suppose $\bar{a}\dind_A \bar{b}$, that is, $\delta(\bar{a}/\bar{b}A)=\delta(\bar{a}/A)$. 
By Lemma~\ref{fiber}, using ($\DC_L$) and ($\FMV$), 
$$\delta(\bar{b}/A)+\delta(\bar{a}/\bar{b}A)=\delta(\bar{a}/A)+\delta(\bar{b}/\bar{a}A).$$
It follows that $\delta(\bar{b}/\bar{a}A)=\delta(\bar{b}/A)$, as required.

\noindent \emph{Algebraic closure}.\enspace Suppose $\chi(\bar{a},\bar{b})$ holds, where $\bar{b}$ is possibly an imaginary tuple in $\acl^{\eq}(A)$, and let 
$\bar{b}=\bar{b}_1,\ldots,\bar{b}_k$ be the conjugates of $\bar{b}$ over $A$. Then $\delta(\chi(\bar{x}/\bar{b}))=\delta(\chi(\bar{x}/\bar{b}_i))$ for each $i$ by ($\DC_L$). There is a
 formula $\rho(\bar{x})\in\tp(\bar{a}/A)$ which is equivalent to $\bigvee_{i=1}^k \chi(\bar{x},\bar{b}_i)$, and by Lemma~\ref{basics}(ii) we
 have $\delta(\rho(\bar{x}))=\delta(\chi(\bar{x},\bar{b}))$.
\edem

\begin{remark} \rm
In Proposition~\ref{d-properties+}(iv), if we assume ($SA$), then the subset $A\subseteq B$ can be taken to be finite: there is a single formula
$\phi(\bar{x},\bar{b})\in \tp(\bar{a}/B)$ such that $\delta(\bar{a}/B)=\delta(\phi(\bar{x},\bar{b}))$, and we may take $A$ to be the set of elements of $\bar{b}$.
\end{remark}

\subsection{Simplicity and forking}\label{simplicityforking}

Here we prove our main results, Theorems~\ref{thm:lowsimple}, \ref{thm:supersimple}, and \ref{thm:forking}, linking ($A$), ($SA$), and 
($\DC_L$) to simplicity and forking, and related notions; these statements are included, respectively, in Theorems~\ref{low2}, \ref{forking2}, and~\ref{supersimple2}.  
Examples are given in Section~\ref{counterexamples} showing that natural strengthenings of the main theorems do not hold.

We first fix some terminology for simple theories, taken from \cite{wagner}.
\bd \ 
Let $T$ be a complete theory, and $M$ an $\omega_1$-saturated model of $T$ from which the parameters below are taken.
\benum
\item A formula $\phi(\ov{x},\ov{y})$ has the \emph{tree property} (with respect to $T$) if there are $k<\omega$ and a sequence 
$\langle \ov{a}_\mu:\mu\in \ ^{<\omega}\omega \rangle$ such that:

(a)  for every $\mu\in\ ^{<\omega}\omega$, the set $\{\phi(\ov{x},\ov{a}_{\mu\smallfrown i}):i<\omega\}$ is $k$-inconsistent, and 

(b) for every $\sigma\in\ ^{\omega}\omega$, the set $\{\phi(\ov{x},\ov{a}_{\sigma\upharpoonright i}):i<\omega\}$ is consistent.

\item The theory $T$ is \emph{simple} if no formula $\phi$  has the tree property with respect to $T$.
\item A \emph{dividing chain of length $\alpha$ for $\phi$}, or {\em dividing $\phi$-chain of length $\alpha$}, is a sequence $(\bar{a}_i:i\in \alpha)$ such that $\bigcup_{i<\alpha}\phi(x,\bar{a}_i)$ is consistent 
and $\phi(\bar{x},\bar{a}_i)$ divides over $\{\bar{a}_j:j<i\}$ for all $i<\alpha$.
\item A simple theory $T$ is \emph{low} if for every formula $\phi$ there is $n_\phi<\omega$ such that there is no dividing $\phi$-chain of length $n_\phi$.
\eenum
\ed

\begin{theorem} \label{low2}
(i) Assume that ($A$) holds. Then $T$ is simple and low.

(ii) If ($A$) and ($\DC_L$) hold then ($\FMV$) holds. 
\end{theorem}

We first note the following lemma, based on  Proposition~\ref{kintersections} and ultimately on Inclusion-Exclusion. 

 
 \begin{lemma} \label{kincon}
 Let $D$ be a $A$-definable subset of $M^r$ in $L$, and let $\phi(\bar{x},\bar{y})$ be an $L$-formula with $l(\bar{x})=r$ and $l(\bar{y})=s$. Let $(\bar{a}_i:i\in I)$ be an infinite $L^+$-indiscernible sequence  over $A$ of elements of $M^s$. Put 
$D_i:= \phi(M^r,\bar{a}_i)$ for each $i\in I$, and suppose that $D_i\subset D$ and $(D_i:i \in I)$ is inconsistent. Then there is some $i\in I$ such that $\delta(D_i)<\delta(D)$.  
 \end{lemma}
 
\bdem  Suppose for a contradiction that $\delta(D_i)=\delta(D)$ for every $i\in I$. Since $I$ is infinite, we may assume without loss of generality that $I=\omega$.

By indiscernibility there is $k\in \omega$ such that $(D_i:i\in\omega)$ is $k$-inconsistent. By indiscernibility in $L^+$ and our assumption that $\delta(D_i)=\delta(D)$ for all $i$, there is some non-zero $m\in\N$ such that $|D_i|>\frac{|D|}{m}$ for each $i$. To ensure that Proposition~\ref{kintersections} is applicable below, we may assume $m\geq 2$.
 
 Consider now the measure $\mu_D$ on $M$, as defined in Section~\ref{firstpreliminaries}. For every $i<\omega$ we have that $\mu_D(D_i)\geq \frac{1}{m}$, and, by Proposition \ref{kintersections}, there are $i_1,\ldots, i_k$ such that
 \[\mu_D\left(\bigcap_{j=1}^k D_{i_j}\right)\geq \dfrac{1}{m^{3^{k-1}}}>0,\]contradicting  $k$-inconsistency.
\edem


\begin{remark} \label{kintersections2}\rm
There is a strengthening of Lemma~\ref{kincon}, in which the assumption that $(D_i:i\in I)$ is inconsistent is weakened to an assumption that for some $k\in\omega$ and all $i_1<\ldots<i_k\in I$, we have $\delta(D_{i_1}\cap \ldots \cap D_{i_k})<\delta(D)$. This is proved by induction on $k$; we omit the details.
\end{remark}

\emph{Proof of Theorem \ref{low2}.} (i)  We show for every  $L$-formula $\phi(\bar{x},\bar{y})$ that every dividing $\phi$-chain has length at most $m:=m_\phi$, where $m_\phi$ is provided by Lemma~\ref{Acontent}; this suffices, by \cite[Proposition 2.8.6]{wagner}.
Suppose for a contradiction that there is a sequence $(\bar{a}_j: 1\leq j\leq m+1)$ from $M^s$
such that each $\phi(\bar{x},\bar{a}_{j})$ divides over $\{\bar{a}_i:i<j\}$. 
 We show by induction that there is a sequence $(\bar{b}_j: 1\leq j\leq m+1)$ such that
$\tp_L(\bar{b}_j: 1\leq j\leq m+1)=\tp_L(\bar{a}_j: 1\leq j\leq m+1)$ and $\delta(\bigwedge_{1\leq i\leq k+1} \phi(x,\bar{b}_i))<\delta(\bigwedge_{1\leq i\leq k} \phi(x,\bar{b}_i))$  for each $k=1,\ldots,m$, contradicting the choice of $m$.

To start the induction, put $\bar{b}_1=\bar{a}_1$. For the induction step, suppose that $\bar{b}_1,\ldots,\bar{b}_k$ have been constructed as above. As
$\tp_L(\bar{b}_j: 1\leq j\leq k)=\tp_L(\bar{a}_j: 1\leq j\leq k)$, there is $\bar{c}$ such that
$\tp_L(\bar{b}_1,\ldots,\bar{b}_k,\bar{c}))=\tp_L(\bar{a}_1,\ldots,\bar{a}_{k+1})$.
We apply Lemma~\ref{kincon} with $A$ the union of the elements in $\{\bar{b}_i:i\leq k\}$ and $D$ the set defined by
$\phi(\bar{x},\bar{b}_1)\wedge \ldots \wedge \phi(\bar{x},\bar{b}_k)$. Let $(\bar{d}_i:i\in \omega)$ be an indiscernible sequence over $A$ 
witnessing the dividing of $\phi(\bar{x},\bar{c})$ over $A$; that is, $\tp(\bar{d}_i/A)=\tp(\bar{c}/A)$ for each $i$ and $\{\phi(\bar{x},\bar{d}_i):i\in \omega)$ is inconsistent. By Ramsey's Theorem, compactness, and $\omega_1$-saturation we may suppose that the tuples $\bar{d}_i$ all lie in $M$ and that the indiscernibility is with respect to $L^+$. Let $D_i$ be the solution set in $M^r$ of $\phi(\bar{x},\bar{b}_1)\wedge \ldots \wedge \phi(\bar{x},\bar{b}_k)\wedge \phi(\bar{x},\bar{d}_i)$ for each $i\in \omega$. By Lemma~\ref{kincon}, there is $i\in \omega$ such that $\delta(D_i)<\delta(D)$. Then put $\bar{b}_{k+1}:=\bar{d}_i$.\hspace{4cm}

(ii) This follows from Lemma~\ref{SOP} and the fact that simplicity implies there is no formula with strict order property. $\Box$

\begin{theorem} \label{forking2}
Let $A, B$ be countable subsets of $M$, and $\bar{a}$ a tuple from $M$. 

(i) Assume ($A$) and ($\DC_L$). Then $$\bar{a}\ind_ AB\Rightarrow \bar{a}\dind_A B.$$

(ii) Assume ($SA$) and ($\DC_L$). Then 
$$\bar{a}\dind_A B\Rightarrow \bar{a}\ind_A B.$$

In particular, under ($SA$) and ($\DC_L$) we have $$\bar{a}\ind_A B\Leftrightarrow \bar{a}\dind_A B.$$
\end{theorem}

We emphasize that this is a statement about forking in the particular ($\omega_1$-saturated) model $M$.  The proof of (ii) uses the following result, which is close to Lemma 2.9 of \cite{h1}. 

\bl \label{forkingdelta} Let $D=\psi(M^r,\ov{a})$ be a definable subset of $M^r$ and $\phi(\ov{x},\ov{b})$ a formula implying $\psi(\ov{x},\ov{a})$. If $\phi(\ov{x},\ov{b})$ divides over $\ov{a}$, then there exists a tuple $\ov{b}'\in M^r$ with $\bar{b}'\models \tp(\ov{b}/\ov{a})$ such that $\delta(\phi(\ov{x},\ov{b}'))<\delta(D)$.
\el

 

\bdem
Since $\phi(\ov{x},\ov{b})$ divides over $\ov{a}$, there is an indiscernible sequence $(\ov{b}_i:i<\omega)$ such that $\ov{b}_i\models \tp(\ov{b}/\ov{a})$ and the set $\{\phi(\ov{x},\ov{b}_i):i<\omega\}$ is $k$-inconsistent for some $k<\omega$. By saturation, we may suppose that this sequence is $L^+$-indiscernible. Since $\ov{b}_i\models \tp(\ov{b}/\ov{a})$, we have $\phi(\ov{x},\ov{b}_i)\subseteq D$. Lemma~\ref{kincon} then yields the result. 
\edem

\emph{Proof of Theorem \ref{forking2}.} (i) We follow the proof of Claim~1 in \cite[Theorem 4.2]{kp2}, which roughly speaking states that any independence relation satisfying the properties from Propositions \ref{d-properties} and \ref{d-properties+} is implied by forking independence. By our assumptions ($\DC_L$) and (A), the properties in these propositions all hold for $\dind$. Note that by ($A$) and Theorem~\ref{low2}, $T$ is simple, so by Lemma~\ref{SOP} we have ($\FMV$).
 
Suppose $\bar{a}\ndind_A B$. Hence $\delta(\bar{a}/AB)<\delta(\bar{a}/A)$
so there is a formula $\phi(\bar{x},\bar{b}) \in \tp(\bar{a}/A B)$ such that $\delta(\phi(\bar{x},\bar{b}))< \delta(\bar{a}/A)$, so $\bar{a}\ndind_A \bar{b}$. We must show
$\bar{a}\nind_A B$, for which it suffices to show $\bar{a}\nind_A \bar{b}$, in particular that $\phi(\bar{x},\bar{b})$ forks over $A$. Suppose for a contradiction that $\bar{a}\ind_A \bar{b}$. By Proposition~\ref{d-properties}(i) and $\omega_1$-saturation, there is in $M$ a sequence $(\bar{b}_i:i\in \omega_1)$ of
 realizations of $\tp(\bar{b}/A)$
such that  $\bar{b}_{i}\dind_A \{\bar{b}_j:j < i\}$ for all $i$.
 By $\omega_1$-saturation and Ramsey's Theorem, using $\DC_L$,
we may suppose $(\bar{b}_i:i\in \omega_1)$ is $A$-indiscernible in $L^+$. 
Let $p=p(\bar{x},\bar{b}):=\tp(\bar{a}/\bar{b}A)$. By our assumption that $\bar{a}\ind_A\bar{b}$, the set
$\bigcup\{\phi(\bar{x},\bar{b}_i):i\in \omega_1\}$ is consistent, realized by $\bar{a}'$, say.
It follows from ($\DC_L$) that $\bar{a}'\ndind_A \bar{b}_i$ for all $i$; indeed, by ($\DC_L$) and ($\FMV$) we have 
$\delta(\phi(\bar{x},\bar{b}_i))=\delta(\phi(\bar{x},\bar{b}))$ for each $i\in\omega_1$. 
Hence,  using symmetry and transitivity 
of $\dind$, it follows that 
$\bar{a}'\ndind_{A\{\bar{b}_j:j< i\}} \bar{b}_i$ for all $i\in \omega_1$: 
indeed, otherwise $\bar{a}'\dind_{A\{\bar{b}_j:j<i\}} \bar{b}_i$, so 
$\bar{b}_i\dind_{A\{\bar{b}_j:j<i\}}\bar{a}'$ and $\bar{b}_i\dind_A \{\bar{b}_j:j<i\}$, so $\bar{b}_i\dind_A \bar{a}'$ and
 hence $\bar{a}'\dind_A \bar{b}_i$, a contradiction.  However, by local character 
of $\dind$, which holds by ($A$), there is $i\in \omega_1$ such that 
$\bar{a}' \dind_{A\{\bar{b}_j:j<i\}} \{\bar{b}_j:j\in \omega_1\}$.
 This is a contradiction.

(ii)
Suppose $\bar{a} \nind_A B$. We must show $\delta(\bar{a}/BA)<\delta(\bar{a}/A)$. By ($SA$), there is $\psi(\bar{x})\in p:=\tp(\bar{a}/A)$ such that $\delta(\psi(\bar{x}))=\delta(p)$. Since $\bar{a}\nind_A B$, there is a formula $\phi(\ov{x},\ov{b})\in \tp(\ov{a}/BA)$ that forks over $A$, and we may suppose that $\phi(M^r,\bar{b})\subseteq \psi(M^r)$. Since (SA) implies (A), by Theorem~\ref{low2} $T$ is simple, so forking and dividing agree. Hence, by 
Lemma~\ref{forkingdelta}, there is some $\ov{b}'\models \tp(\ov{b}/A)$ such that $\delta(\phi(\ov{x},\ov{b}'))<\delta(\psi(\ov{x}))$. If $\alpha\in \Aut(M/A)$ is such that $\alpha(\ov{b})=\ov{b}'$, then
$\phi(\alpha(\bar{a}),\bar{b}')$ holds. Hence $\alpha(\ov{a})\nind^\delta_A \ov{b}'$, which implies by invariance of $\dind$ (using ($\DC_L$)) that $\ov{a}\nind^\delta_A \ov{b}$. 
\qedsymbol

\begin{remark} \label{SAnotOD}\rm
1. The above proof shows, assuming ($SA$) and ($\DC_L$), that a formula $\phi(\bar{x},\bar{b})$ forks over $A$ if and only if, for every $\psi(\bar{x})\in L_A$ that is consistent with $\phi$, we have $\delta(\phi\wedge\psi)<\delta(\psi)$. The direction 
($\Leftarrow$) again just requires ($A$) and ($\DC_L$).

2. The proof of (ii) yields the following, just under the assumption ($SA$), so without ($\DC_L$). Assume ($SA$) and suppose 
$\bar{a} \nind_A B$. Then there is $\bar{a}'B'$ such that $\tp_L(\bar{a}B/A)=\tp_L(\bar{a}'B'/A)$ and
$\bar{a}'\ndind_A B'$.
Indeed, in the  proof above and with $\phi(\bar{x},\bar{b})$ as above, we obtain $\bar{b}'$ such that
$\delta(\psi(\bar{x})\wedge \phi(\bar{x},\bar{b}'))<\delta(\psi(\bar{x}))$. Choose  $\bar{a}'B'$ by $\omega_1$-saturation so that
$\tp_L(\bar{a}\bar{b}B/A)=\tp_L(\bar{a}'\bar{b}'B'/A)$. Then $\bar{a}' \nind^\delta_A B'$.
The point to note here is that $\delta$-dimension is not part of the $L$-type, and is not preserved in general by automorphisms of the ultraproduct.

3. Example~\ref{forkingSA} shows that (ii) is not true if either of the assumptions $(SA)$ or $(\DC_L)$ is dropped.
\end{remark}

\begin{question} \rm
Is there a local version of Theorem~\ref{forking2}? For example, in the setting ($\DC_L$) and ($\FMV$), along with $A_\phi^*$, is it true that for all $\bar{a},\bar{b}$ and countable $A$,
some Boolean combination $\psi(\bar{x})$ of $\phi$-formulas in $\tp(\bar{a}/A\bar{b})$ forks over $C$ if and only if 
$\delta(\psi(\bar{x}))<\delta(\chi(\bar{x}))$ for every Boolean combination $\chi$ of $\phi$-formulas in $\tp(\bar{a}/A)$?
\end{question}

Most parts of the proof of Theorem~\ref{forking2} localize easily, to Boolean combinations of instances of a finite set $\Gamma$ of formulas $\phi(\bar{x},\bar{y})$. The problem occurs in the proof that $\delta$-forking implies forking, when symmetry of $\dind$ is applied. This rests on Lemma~\ref{fiber}, which we have not been able to localize, essentially because of the additional formulas used to witness ($\DC_L$).

\begin{theorem} \label{supersimple2}
Assume ($SA$).
Then $T$ is supersimple.
\end{theorem}

\bdem Suppose for a contradiction that $T$ is not supersimple. Then there are countable sets $B_0\subset B_1\subset \cdots$ and a type $p$ over $\displaystyle{B=\bigcup_{i<\omega}B_i}$ such that for all $i<\omega$, $p|_{B_{i+1}}$ forks over $B_i$. Let $\bar{a}$ realise $p$.

For every $n$ we build sets $B_n'$, with $B_0'\subset\ldots\subset B_n'$, along with tuples $\bar{a}_n$, such that
$\tp_L(\bar{a}_n B_n')=\tp_L(\bar{a} B_n)$  and $\delta(\bar{a}_n/B_{i+1}')<\delta(\bar{a}_n/B_i')$ for each $i<n$. Suppose that these have been found for some given $n$. Then there is $B_{n+1}^*$ such that $\tp_L(\bar{a}_n/B_0'\ldots B_n'B_{n+1}^*)=\tp_L(\bar{a}B_0\ldots B_nB_{n+1})$. Thus $\bar{a}_n\ind_{B_0'\ldots B_n'} B_{n+1}^*$, so by Remark~\ref{SAnotOD}(2), there is
$\bar{a}_{n+1}B_{n+1}'$ such that $$\tp_L(\bar{a}_{n+1}B_{n+1}'/B_0'\ldots B_n')=\tp_L(\bar{a}_nB_{n+1}^*/B_0'\ldots B_n')$$ and
$\delta(\bar{a}_{n+1}/B_{n+1}')<\delta(\bar{a}_{n+1}/B_0'\ldots B_n')$. 

Finally, let $B':=\bigcup(B_i':i\in \omega)$ and put $p_n:=\tp(\bar{a}_n/B_n')$ for each $n$, and $p':= \bigcup(p_n:n\in \omega)$. Then $p\in S(B')$ and $\delta(p'|B_{n+1}')<\delta(p'|B_n')$ for each $n$, contradicting $(SA)$. (In fact, this contradicts the weaker assumption $(SA^-)$, but we need ($SA$) when using Remark~\ref{SAnotOD}(2); see also Example~\ref{forkingSA}.) \edem

\subsection{Pseudofinite dimension and stability}

In Proposition~\ref{stablecrit} below we characterize,  among structures $M$ satisfying ($A$) --- so among simple structures $M$---  when $M$ is stable. 
The statement involves the measure $\mu_D$ on definable subsets of $D$ defined in Section~\ref{firstpreliminaries} as well as dimension, and is presented as a local result. Note that by Lemma~\ref{A and A*}, the assumption ($A_\phi^*$) for all $\phi$ is equivalent to ($A$). In Example \ref{stablenonattainability}
we present an example demonstrating that stability  does not necessarily imply $(A)$.

\begin{proposition}\label{stablecrit}
Assume $(A)_\phi^*$. Then the following are equivalent:

{\rm (i)} $\phi(\bar{x},\bar{y})$ is unstable;

{\rm (ii)} there is for some $\bar{d}$ a $\bar{d}$-definable set $D\subseteq M^r$ and a sequence $(\bar{a}_i:i\in\omega)$, $L^+\/$-indiscernible over 
$\bar{d}$,  such that
$\delta(D)=\delta (D\wedge \bigwedge_{i\in\omega}\phi(\bar{x},\bar{a}_i))$,  and 
$$\mu_D(\phi(\bar{x},\bar{a}_i)\wedge \phi(\bar{x},\bar{a}_j))<\mu_D(\phi(\bar{x},\bar{a}_i))\mbox{\ for all $i<j$.}$$
\end{proposition}

\bdem (i) $\Rightarrow$ (ii)\enspace Assume~(i), and find a  set $D$, defined by a finite partial $\phi$-type, such that there are 
$\{\bar{a}_i:i\in\omega\}$ and $\{\bar{b}_i:i\in\omega\}\subset D$ with $\phi(\bar{b}_i,\bar{a}_j)$ holding if and only if $i>j$, and such that there is no such set $D'\subset D$ with $\delta(D')<\delta(D)$. This is possible by $(A)_\phi^*$. 

Suppose that (ii) is false, and that $(\bar{a}_i\bar{b}_i:i\in\omega+1)$ is $L^+$-indiscernible over parameters used to define 
$D$,  with $\bar{b}_i\in D$ for all $i\in \omega+1$, and $\phi(\bar{b}_i,\bar{a}_j)$ holding if and only if $i>j$. It follows from indiscernibility that $\mu_D(\phi(\bar{x},\bar{a}_i))$ is constant. 
By the minimality in the choice of $D$, we have $\delta(D)=\delta(D\cap \phi(M^r,\bar{a}_0))$. It follows again by $L^+$-indiscernibility that $\delta(D)=\delta(D\cap \phi(M^r,\bar{a}_i))$ for each $i$; for there is some fixed $N\in \omega$ such that $|D|<N|D\cap \phi(M^r,\bar{a}_i)|$.

As we have assumed that~(ii) is false, by $L^+$-indiscernibility for all $i<j < \omega_1$ we have 
$\mu_D(\phi(\bar{x},\bar{a}_i)\wedge \phi(\bar{x},\bar{a}_j))=\mu_D(\phi(\bar{x},\bar{a}_i))$. 
Hence  
$\mu_D(\phi(\bar{x},\bar{a}_i)\wedge \neg \phi(\bar{x},\bar{a}_j))=0$, so $\delta(D \cap \phi(\bar{x},\bar{a}_i)\wedge \neg \phi(\bar{x},\bar{a}_j))<\delta(D)$. 
 Now put 
$D'=D\cap(\phi(\bar{x},\bar{a}_0)\wedge \neg\phi(\bar{x},\bar{a}_\omega))$. Then $\delta(D')<\delta(D)$. Since $\bar{b}_i\in D'$ for each $i>0$, this contradicts the assumption of minimality of $\delta(D)$ in the choice of $D$.  

(ii) $\Rightarrow$ (i)\enspace We employ the same strategy used in the proof of (ii) $\Rightarrow$ (i) in \cite[Proposition 4.2]{ms}.
Assume that~(ii) holds, and for each $i$ let $S_i=D\cap \phi(\bar{x},\bar{a}_i)$.  By indiscernibility and 
Proposition~\ref{kintersections}, for all $k$ there is $\epsilon_k$ with
$0<\epsilon_k<1$ such that for all $i_1<\ldots<i_k\in\omega$, we have $\mu_D(S_{i_1}\cap\ldots\cap S_{i_k})=\epsilon_k$. By assumption, $\epsilon_2<\epsilon_1$, hence for $i<j$ the set $S_i\setminus S_j$ has positive measure $\epsilon_1-\epsilon_2$.

We claim that for all $1\leq p\leq q\in\omega$, there is a positive 
real number $\rho=\rho(p,q)$
such  that
\[
(*)\hspace{1cm}\mu_D\left(\left\{\bar{x}: \bigwedge_{1\leq i\leq p}\phi(\bar{x},\bar{a}_i)
\wedge
\bigwedge_{p+1\leq j\leq q}
\neg\phi(\bar{x},\bar{a}_{j})\right\}\right) =\rho.
\]

We prove this by induction on $d=q-p$. Observe that~($*$) is true for the pairs $(1,1)$, more generally, for all pairs 
$(p,p)$, and for $(1,2)$, with $\rho(p,p)=\epsilon_p$ and $\rho(1,2)=\epsilon_1-\epsilon_2$. 

We first show by induction on $p$ that~($*$) holds for all pairs $(p,p+1)$ with $p\geq 1$. For the induction step, let $N=N(p,p+1)$ be sufficiently large (a number  guaranteed to exist by 
Remark~\ref{k=2} applied to the sets $T_j$ below). For $p\leq j\leq p+N$ put 
\[
T_j:=\bigwedge_{1\leq i\leq p-1} \phi(\ov{x},\ov{a}_i)\wedge   \phi(\ov{x},\ov{a}_j)  \wedge \neg \phi(\ov{x};\ov{a}_{p+N+1}).  
\] 
By induction hypothesis and indiscernibility, $\mu_D(T_j)=\rho(p,p+1)$ for all $j$ with $p\leq j\leq p+N$. Hence, by the choice of $N$, we have $\mu_D(T_j\cap T_{j'})>0$ for all $j,j'$ with $p\leq j,j'\leq p+N$.  Thus, by indiscernibility, ($*$) holds for all pairs $(p+1,p+2)$, completing this first induction. 

Now suppose that~($*$) is true
for all pairs $(p,q)$ with $q-p=d\geq 1$. We prove that it holds for all pairs $(p, q+1)$. 
For $k\in\omega$ with $q\leq k$, put 
\[
U_k:=\left(\bigwedge_{1\leq i\leq p} \phi(\ov{x},\ov{a}_i)\>\wedge \bigwedge_{p+1\leq j\leq q-1} \neg \phi(\ov{x};\ov{a}_j)\right) \wedge \neg \phi(\ov{x};\ov{a}_k). 
\] 
By $L^+$-indiscernibility and the induction hypothesis, there is a number $\rho=\rho(p,q)$ such that 
$\mu_D(U_k)=\rho$ for all $ k\in\omega$ with $q\leq k$.
By Proposition~\ref{kintersections} once again there are $k,k'$ with $q\leq k<k'$ such that $\mu_D(U_k\cap U_{k'})\geq \rho^3>0$. By indiscernibility, we may take $\rho(p,q+1)$ to be $\mu_D(U_k\cap U_{k'})$, which completes the induction and thus the proof of~($*$).


Given~($*$), as $\rho>0$, it is easy to find $(\bar{b}_i:i\in\omega\setminus\{0\})$ such that $\phi(\bar{b}_i,\bar{a}_j)$ holds if and only if $i>j$.\edem

\begin{remark} \rm 
Proposition~\ref{asymp} below shows that ultraproducts of asymptotic classes (in the sense of \cite{ms} and \cite{elwes}) provide natural examples satisfying ($SA$) and ($\DC_L$). It was shown in  \cite[Proposition 6.5]{elwes} (with a point of confusion concerning unimodularity clarified in \cite{kestner}) that any {\em stable} ultraproduct of an asymptotic class --- in fact, any stable structure which is measurable in the sense of \cite{ms} --- is one-based. It would be interesting to generalize this to the assumptions of this paper.

It is not true that every pseudofinite superstable structure is one-based. For example, in \cite[Section 5]{mt0}, an example is sketched of a pseudofinite $\omega$-stable group $G$ that is not nilpotent-by-finite. It has the form $({\mathbb C},+,\cdot,T)$ where $T$ is an infinite subgroup of $({\mathbb C}^*,\cdot)$, and its construction is due, independently, to Chapuis, Khelif, Simonetta, and Zilber. By the main theorem of \cite{hp0}, one-based stable groups are abelian-by-finite, so $G$ is not one-based. We have not checked if it satisfies conditions such as ($A$) or ($\DC_L$). There is also an example in \cite[Section 5]{mt0} of a nilpotent class 2 but not abelian-by-finite $\omega$-stable pseudofinite group $G$, based on the Mekler construction that codes graphs into groups; again $G$ cannot be one-based.
\end{remark}

\section{Examples}\label{examples}

This section has two aims. We first present examples designed to show that the obvious strengthenings of the principal results in Section~\ref{simplicityforking} fail. Then in \ref{asymptotic}, \ref{pseudofinitevs}, and \ref{homocyclic} we investigate the conditions ($A$), ($SA$), and ($\DC_L$) in the context of rather natural examples: asymptotic classes of finite structures; pseudofinite  2-sorted infinite-dimensional vector spaces over pseudofinite fields; and, ultraproducts of homocyclic groups. 

\subsection{Counterexamples}\label{counterexamples}


\bej \label{convsupersimple} \rm We show that the converse to Theorem~\ref{supersimple2}  is false.  Let $L$ consist of infinitely many unary predicates $(P_i:i\in \omega)$. For each $n$, let $M_n$ be a finite structure with domain of size $n^n$ 
such that $P_i(M_n)\supseteq P_{i+1}(M_n)$ for all $i$, and $|P_i(M_n)|=n^{n-i}$ for $i\leq n$ and $P_i(M_n)=\emptyset$ for 
$i>n$. Let $M$ be a non-principal ultraproduct. Then $\Th(M)$ is superstable of U-rank 1. However, 
$\delta(P_i(M))>\delta(P_{i+1}(M))$ for all $i$, so ($SA$) fails.

We further note  that supersimplicity does not follow from ($A$) + ($\DC_L$): Proposition~\ref{abelian} below provides a counterexample. Also, Proposition~\ref{329converse} provides a partial converse to Theorem~\ref{supersimple2} for expansions of groups.
\eej

\begin{problem} \rm.
Find natural conditions on pseudofinite dimension in $M$ that are equivalent to supersimplicity of $\Th(M)$.
\end{problem}

\begin{example} \label{stablenonattainability} From the definitions one might believe that stability of a formula $\phi(x,y)$ implies $(A_\phi)$, but we present here a counterexample showing that not even the stability of $\Th(M)$ implies attainability. 

For each $n<\omega$, let $(M_n,E)$ be the structure where $M_n$ has $\displaystyle{\sum_{i=1}^n n^i}$ elements and $E$ is interpreted as an equivalence relation with a class of size $n^i$ for each $i\leq n$. Let $a_{n,i}\in M_n$ be an element in the class of size $n^{n-i}$. The theory of $M=\prod_{\mathcal{U}}M_n$ is just the theory of an equivalence relation with infinitely many infinite classes, so $\Th(M)$ is stable of $U$-rank $2$ (and low). By taking the formula $\phi(x,y):=\neg (xEy)$ and the sequence of elements $\langle a_i:=[(a_{n,i})_{n<\omega}]_{\mathcal{U}}:i<\omega\rangle$ in $M$, we show that attainability fails for the positive $\phi$-type $p:=\{\phi(x,a_i):i<\omega\}$. 

First note for every $t<\omega$ that $\delta(xEa_t)>\delta(xEa_{t+1})$. Otherwise, there would be a natural number $N$ such that, for $\mathcal{U}$-almost all $n$, we have: 
\begin{eqnarray*}
\log ( |xEa_{n,t}| ) - \log (|xEa_{n,t+1}|)\leq N;\\
\log (n^{n-t})-\log (n^{n-(t+1)})\leq N; \\
(n-t)\cdot \log n - (n-t-1)\cdot \log\,n\leq N;\\
\mbox{and hence\ }\log n\leq N,
\end{eqnarray*} 
which contradicts the fact that $\{n<\omega: n\leq N\}\not\in\mathcal{U}$. A similar argument shows for every $t<\omega$ 
that $\delta(xEa_t)>\delta(p)$. 

Therefore, given finitely many $a_{i_1},\ldots,a_{i_k}$ the set defined by $\bigwedge_{j\leq k} \phi(x;a_{i_j})$ contains $xEa_{t}$ for $t>i_1,\ldots,i_k$, and we have\[\delta\left(\bigwedge_{j=1}^k \phi(x;a_{i_j})\right)\geq \delta(xEa_t)>\delta(xEa_{t+1})>\delta(p).\] Thus, ($A_\phi$) does not hold.

This example also shows that the converse to Theorem~\ref{low2}(i) fails.
\end{example}


The next example is close in spirit to that discussed in Theorem~\ref{cdmvs}, and will be developed more fully in subsequent work. 

\begin{example} \rm \label{bilinear}
Fix a prime $p$, and let $\C$ be the collection of all 2-sorted structures consisting of a finite field $F$ of characteristic $p$ and an even-dimensional vector space $V$ over $F$, with a function 
$F\times V\to V$ for scalar multiplication and also a function symbol $\beta:V\times V\to F$ interpreted by a non-degenerate alternating bilinear form, that is, a symplectic form.  Let $M=(V^*,F^*)$ be an ultraproduct of members of $\C$ with both $F^*$ and $\dim V^*$ infinite. By Proposition 7.4.1 of \cite{granger},
$\Th(V^*,F^*)$ is not simple. 

It is easy to see directly that $M$ does not satisfy ($A$). Indeed, let $\phi(x,y)$ be the formula $\beta(x,y)=0$. For $(V,F)\in \C$ and $a_1,\ldots,a_n\in V$ that are linearly independent with $\beta(a_i,a_j)=0$ for all $i,j\in \{1,\ldots,n\}$, let $X_n:=\{x: \phi(x,a_1)\wedge \ldots\wedge\phi(x,a_n)\}$. Then
$|X_n|=\frac{|X_{n-1}|}{|F|}$. Thus, if sets $X_n'$ are defined in the same way in the ultraproduct, we would have $\delta(X_n)>\delta(X_{n+1})$ for all $n$. 

In work in preparation by the second and third authors with W. Anscombe and D. Wood, it is shown that $M$ has ($\DC_L$); in fact, an asymptotic result analogous to Theorem~\ref{cdmvs} below is proved. It can be seen that if $a\in V^*\setminus\{0\}$ then the formula $\beta(x,a)=0$ does not fork over $\emptyset$. Thus, in Theorem~\ref{forking2}(i), the condition ($A$) cannot be omitted.
\end{example}
 
\begin{example}\label{forkingSA} \rm We now show that the assumption ($SA$) cannot be omitted in Theorem~\ref{forking2}(ii), or weakened to $(SA^-)$. Let $L$ be a language with binary relations $(E_i)_{i \in \omega}$,  all to be interpreted by equivalence relations, a binary relation $F$, and a unary relation $P$. In the ultraproduct $M$ of a class $\mathcal C=\{M_j: j<\omega\}$ of finite $L\/$ structures we want for each $i<\omega$ that each $E_i$-class is a union of infinitely many $E_{i+1}$-classes. Let $E$ be the intersection of the relations $E_i$ in $M$, a 
$\bigwedge$-definable equivalence relation. The predicate $P$ in $M$ is to have infinite coinfinite intersection with each $E$-class. Lastly, the interpretation of $F$ is an equivalence relation on the complement of $P$ with infinitely many classes,  such that  each $F$-class has infinite intersection with each $E_i$-class for each $i<\omega$. 
We can arrange that the finite $L$-structures $M_j$ are built in a uniform way, so that, for example, if $j$ is sufficiently large relative to $i$ then each $E_i$-class of $M_j$ is a union of $E_{i+1}$-classes all of the same size, with a corresponding uniformity for $P$ and $F$. We also can ensure for each $i$ that all $E_i$-classes have the same pseudofinite dimension $\delta_i$, and that if $Q$ is an $E_i$-class and $D$ an $F$-class then $\delta_i>\delta(Q\cap D)>\delta_{i+1}$. 

The theory of $M$ is stable, has quantifier elimination, but is not superstable, so $(SA)$ fails. The finite structures can be chosen so that $M$ satisfies ($A$) and ($\DC_L$), and also $(SA^-)$. Let $c$ realize $P$, let $q$ be the type over $\{c\}$ containing $\{\neg P(x), E_i(x,c):i\in \omega\}$, and let $b$ realize $q$. Then clearly there is a type $q'$ over $\{c,b\}$ containing $F(x,b)$ and extending $q$ that forks over $c$ but satisfies $\delta(q')=\delta(q)$.

The condition ($\DC_L$) cannot be dropped in Theorem~\ref{forking2}(ii) either. Indeed, consider the following variant of  Example~\ref{earlyexample} (see also the second example in \ref{findelta} below). The language $L$ consists
of a single binary relation $E$, and the class $\C$  consists of finite structures $M_n$ for $n>0$ in which $E$ is interpreted by an equivalence relation with $n$ classes of size $n$ and one of size $n^2$. The ultraproduct $M$ is superstable of rank 2 and has (SA). There is a unique equivalence class $B$ of $M$ such that $\delta(B)=\delta(M)$. 
If $b\in B$, then the formula $E(x,b)$ forks over $\emptyset$, and if $b'\in B\setminus \{b\}$ then $b'\nind_\emptyset b$ but $b'\dind_\emptyset b$. 
This example shows also that in Lemma~\ref{forkingdelta}, we cannot expect the conclusion $\delta(\phi(\bar{x},\bar{b}))<\delta(D)$. 

\end{example}

\begin{example} \rm  \label{findelta} While Theorem~\ref{low2}(ii) asserts that ($A$) + ($\DC_L$) implies ($FMV$), the conditions ($SA$) + ($\DC_L$) do not imply that the set $\{\delta(X): X\subset M, X\mbox{~definable}\}$ is finite. 
Theorem~\ref{cdmvs}, below, provides an example: using the notation there, in the vector space sort, for each $k$ there is a $k$-dimensional definable subspace $V_k$, and for $k<l$ we have $\delta(V_k)<\delta(V_l)$.

We also have an example of a structure with finite SU-rank that satisfies $SA$ but not ($\DC_L$) in which $\delta$ takes infinitely many values on a  uniformly definable family. Let $L$ be a language with a single binary relation $E$, and for each $n<\omega$ let $M_n$ be an $L$-structure with $n\displaystyle{\sum_{i=1}^n n^i}$ elements, with $E$ interpreted by an equivalence relation  having, for each $i\in\{1,\ldots,n\}$, exactly $n$ equivalence classes of size $n^i$. In the ultraproduct $M$, the equivalence relation $E$ has infinitely many classes, all infinite, so $\Th(M)$ is $\omega$-stable of Morley rank 2. However the {\em uniformly} definable family of sets $\{x: Exb\}$ of equivalence classes takes infinitely many $\delta$-values as $b$ varies. It is routine to check using quantifier elimination that ($SA$) holds, and so ($\DC_L$) fails by Theorem~\ref{low2} (ii).

\begin{example} 
Under the assumptions ($SA$) and ($\DC_L$), we may have a definable set $X$ in $M$ and a definable subset $Y\subseteq X$ such that $\delta(Y)=\delta(X)$ and $SU(Y)<SU(X)$. To see this, let $L$ have a unary predicate $P$ and a binary predicate $E$.  For $n<\omega$ let $M_n$ be a finite structure in which approximately half the elements satisfy $P$, and $E$ is interpreted by an equivalence relation on $P$ with classes all of the same size, with both the size and number of $E$-classes increasing as  $n\to \infty$. If $M$ is a non-principal ultraproduct of $(M_n : n<\omega)$ then  $\delta(\neg P(x))=\delta(M)$ but $SU(\neg P(x))=1<2=SU(M)$. 

We note here that Lemma~\ref{SUdelta} shows that this phenomenon is not possible if $X$ has a definable group structure.
\end{example}

\end{example}

\begin{remark} \rm
The elimination of imaginaries assumption is required in Proposition~\ref{sopei}(ii). 
Consider the class $\C$ of finite structures $M_n$ equipped with binary relations $E$ and $<$ such that $E$ is an equivalence relation with $n$ classes totally ordered by $<$ (so $M_n$ is a preorder). Furthermore, as $n$ increases, the size of all classes should increase without bound, and for $k<k+1\leq n$, the size of the $<$-$k+1$st equivalence class should be much larger than the $n$-fold cartesian product of the $k$th class.  The ultraproduct $M$ of $\C$ is a discretely ordered preorder whose quotient is infinite and has first and last elements---and thus has the strict order property---such that each equivalence class is infinite. Moreover, $M$ admits elimination of quantifiers in a suitable enlargement of the language and evidently has ($\DC_L$).
\end{remark}



\subsection{Asymptotic classes}\label{asymptotic}

As usual, we consider a class $\C$ of finite structures in a language $L$, with a corresponding class $\C^+$ in the 2-sorted language $L^+$ of which we form an ultraproduct $K^{*}$ (often denoted $K^{*}(\C)$) in the 2-sorted language $L^+$, and then consider the reduct $M$ to $L$ of the structure induced by $K^*$ on the sort ${\bf D}$. 

Recall from \cite{elwes} (see also \cite{ms} and \cite{em}) the definition of an {\em asymptotic class} of finite structures, and the corresponding notion of {\em measurable structure}. Every non-principal ultraproduct of an asymptotic class is measurable---but not conversely---and every measurable structure has supersimple finite rank theory. The motivating example of an asymptotic class, by the main theorem of \cite{cdm}, is the class of finite fields. Likewise, by \cite[Theorem 3.5.8]{ryten} of Ryten, if $p$ is a prime and $m,n$ are coprime natural numbers with $m>1$ and $n\geq 1$,
then the class of all finite difference fields $({\mathbb F}_{p^{kn+m}}, {\rm Frob}^k)  $, where ${\rm Frob}$ is the Frobenius automorphism $x\mapsto x^p$, is a 1-dimensional asymptotic class. Using this, Ryten showed that the class of all finite simple groups of any fixed Lie type is an asymptotic class. Likewise, Elwes \cite{elwes} showed that any smoothly approximable structure  has an approximating sequence of envelopes which forms an asymptotic class.

Our first result shows that asymptotic classes and their ultraproducts fall under the framework of Section~\ref{qfd}.
If $\C$ is an $N$-dimensional asymptotic class, $M_n \in \C$, and $X$ is a definable set in $M_n$ of cardinality approximately $\mu|M_n|^d$ (in the sense of asymptotic classes), we shall say that $X$ has dimension $Nd$. This notion is well-defined, provided $M_n$ is sufficiently large, relative to the formula defining $X$.

\begin{proposition} \label{asymp}
Let $\C$ be an asymptotic class of finite structures in a language $L$. Then $K^{*}(\C)$ 
satisfies ($SA$), ($\DC_L$) and ($FMV$),  and for all definable subsets $X, Y\subseteq M^r$, $\dim(X)=\dim(Y)$ if and only if $\delta(X)=\delta(Y)$.
\end{proposition}

\bdem We show first that $K^{*}(\C)$ satisfies ($SA$). So suppose that ($SA$) is false, and let $M$ be the corresponding ultraproduct of $\C$.
There is a countable  set $C\subset M$, some $p\in S(C)$,  and a sequence of formulas $\phi_i(x)\in p$ for $i\in \omega$ such that $\delta(\phi_{i+1}(x))<\delta(\phi_i(x))$ for all $i$. We may suppose in addition that 
$M\models \forall\bar{x}(\phi_{i+1}(\bar{x})\to \phi_i(\bar{x}))$.  In particular, for each $i,n< \omega$, we have
$\frac{|\phi_i(x)|}{|\phi_{i+1}(x)|} >n$ (non-standard cardinalities). In the class $\C$, this means that for each $n$ there is a set $U$ in the ultrafilter such that for all $j\in U$ and $i \leq n$, we have  
$n|\phi_{i+1}(M_j^r)|<  |\phi_{i}(M_j^r)|$. 
Recall (see \cite{em}) that in the ultraproduct $M$ of the asymptotic class $\C\/$, dimension in the sense of measurable structures, which we denote here by $\dim$, is induced by dimension in the asymptotic class. We thus have 
$\dim(\phi_{i+1}(M^r))<\dim(\phi_i(M^r))$ for each $i< \omega$. 
It follows that in the ultraproduct $M$, if $\dim$ denotes the dimension in the sense of measurable structures, we have 
$\dim(\phi_{i+1}(M^r))<\dim(\phi_i(M^r))$ for each $i< \omega$. Since dimensions in an asymptotic class are non-negative integer valued, and the dimension of any subset of $M^{l(x)}$ is at most
$N^{l(x)}$ where $\C$ is an $N$-dimensional asymptotic class, this is a contradiction. 

Next, we prove the final assertion of the proposition, from which ($\DC_L$) follows, as dimension is definable in measurable structures. Let $X,Y\subset M^r$ be definable. Then 
$$\dim(X)=\dim(Y)\Leftrightarrow \exists n<\omega\> \left(\frac{1}{n}\leq \frac{|X|}{|Y|}\leq n\right)
\>\Leftrightarrow\>\delta(X)=\delta(Y),$$
where the first equivalence follows from the behavior of dimension in asymptotic classes by taking $n\geq \frac{1}{\mu}$ where 
$\mu$ is the minimum of the finitely many positive values for the measures of formulas defining $X$ and $Y$.

Finally, by Theorem~\ref{low2}(ii), ($FMV$) follows from ($A$) and ($\DC_L$).
\edem

\begin{remark} \rm
It is known (I.\ Ben-Yaacov, personal communication) that in a measurable structure $M$, change in dimension corresponds to forking. That is, $\dim(a/B\cup C)=\dim(a/C)$ if and only if ${\rm SU}(a/B \cup C)={\rm SU}(a/C)$. However, we do not claim that in every measurable structure  
the dimension can be adjusted so that it coincides with SU-rank. As an example, in a language $L$ with a single unary predicate $P$,
 consider the class $\C$ of finite $L$-structures $M_n$ of size $n^2$ for $n<\omega$ such that $P(M_n)$ has size $n$. This is a 2-dimensional asymptotic class, and taking an ultraproduct, the universe has dimension~2 but the SU-rank is~1. Although the dimension of $x=x$ could be changed to~1, this could not be done preserving the relation `has the same dimension'.
\end{remark}

\begin{question} \rm 
If $\C$ is an asymptotic class, we may consider the class $\C$ of corresponding 2-sorted structures in the language $L^+$, and the ultraproducts $K^*(\C)$. In key examples such as if $\C$ is the class of finite fields, does $K^*(\C)$ have an NTP2 theory?
\end{question}

\subsection{Pseudofinite vector spaces}\label{pseudofinitevs}

In this and the next section, we consider two examples that belong to the framework of `multi-dimensional asymptotic classes' being developed by the second and third authors along with W.~Anscombe and D.~Wood. These are classes of finite structures, typically in a multi-sorted language, in which there is a strong uniformity in the asymptotic cardinalities of definable sets, in terms of the cardinalities of certain sorts. The example in Theorem~\ref{cdmvs} below is prototypical, and, unlike ultraproducts of asymptotic classes, has infinite SU-rank. Example~\ref{bilinear} is similar, but has ultraproducts whose theory is not simple. The example in Proposition~\ref{abelian} also belongs to this framework, but has the additional property that the cardinalities of definable sets are given exactly rather than asymptotically. We also take the opportunity in this section to include some further observations on the model theory of infinite-dimensional vector spaces over pseudofinite fields in Propositions~\ref{nmc} and \ref{genstab}.

Let  $L$ be a  2-sorted language containing: a sort $\mathcal{V}$, the vector space sort, equipped with  a binary function symbol +, a unary function symbol -, and a constant symbol 0; a sort $\mathcal{K}$, the field sort, equipped with the language $L_R$ of rings; and, a function symbol for  scalar multiplication $\mathcal{K} \times \mathcal{V} \to \mathcal{V}$. 
Let $L_{vs}$ be obtained from $L$ by adding, for each $n>0$, an $n$-ary relation symbol $\theta_n$. In a vector space $V$ over $K$ we interpret $\theta_n(v_1,\ldots,v_n)$ as expressing that the vectors $v_1,\ldots,v_n$ are linearly independent. Let $T_{vs}$ be the theory of infinite-dimensional vector spaces, in this language $L_{vs}$, and for a field $F$ let $T_{vs}(F)$ be the theory of infinite-dimensional vector
 spaces over models of $\Th(F)$. It is clear that each theory $T_{vs}(F)$ is complete (see also Granger \cite[Corollary 11.1.6]{granger}). Furthermore, $T_{vs}(F)$ eliminates quantifiers in the  sort $\mathcal{V}$. More formally, by Proposition 11.1.7 of \cite{granger}, which itself is a slight elaboration of the main theorem of Kuzichev \cite{kuzichev} (see also \cite{pierce}), we have the following, where a {\em $\theta$-formula} is an instance of some $\theta_n$.

\begin{lemma} \label{kuz-pierce} Let $\phi$ be a formula of $L_{vs}$. Then $\phi$ is equivalent modulo $T_{vs}(F)$ to a boolean combination of $\theta$-formulas, quantifier-free formulas, and $L_R$-formulas in the field sort.
\end{lemma}

Our main result in this section is the following.

\begin{theorem}\label{cdmvs}
Let $\C$ be the class of all $L_{vs}$-structures $(V,F)$ where $V$ is a finite-dimensional vector space over the finite field $F$.  Let $\phi(x_1,\ldots,x_r,y_1,\ldots,y_s)$ be a formula, and ${\bf V}, {\bf F}$ be indeterminates. Then there is a  finite set $E$ of
polynomials $p({\bf V},{\bf F})\in {\mathbb Q}[{\bf V},{\bf F}]$,  
such that for every $M=(V,F) \in \C$ and all $a_1,\ldots,a_s\in M$, there is some $p({\bf V},{\bf F})\in E$ such that 
\begin{equation*}
\big||\phi(M^r,a_1,\ldots,a_r)| - p(|V|,|F|)\big|= o(p(|V|,|F|) \tag{$*$}
\end{equation*} 
Furthermore---the definability condition---for all $p\in E$, there is a formula $\phi_p(y_1,\ldots,y_s)$
such that if $M\in \C$ is sufficiently large and $a_1,\ldots,a_s\in M$,
then $(*)$ holds if and only if
$M\models \phi_p(a_1,\ldots,a_s)$.
\end{theorem}

The following corollary is immediate, as in Proposition~\ref{asymp}.

\begin{corollary} 
Let $M$ be an infinite ultraproduct of members of $\C$, viewed as usual as the structure induced on the sort ${\bf D}$ (itself formally a pair of sorts) by an ultraproduct $K^*(\C)$. 
Then $M$ satisfies (SA) and $(\DC_L)$.
\end{corollary}

The proof of Theorem~\ref{cdmvs} is by induction on $r$, using a fibering argument in spirit like the proof of the o-minimal Cell Decomposition Theorem. The main work, which rests on \cite{cdm} and hence ultimately on the Lang-Weil estimates, is contained in the following lemma that starts the induction.

\begin{lemma} \label{r=1}
The conclusion of Theorem~\ref{cdmvs} holds when $r=1$, that is, for formulas $\phi(x,y_1,\ldots,y_s)$.
\end{lemma}

\bdem We consider separately the two cases where $x$ lies in the sort $\mathcal{V}$ or the sort $\mathcal{K}$. We work in a large finite structure $M=(V,F)\in \C$.

\smallskip

\noindent {\em Case 1.}\enspace Suppose first that $x$ lies in the sort $\mathcal{V}$, and, replacing $x$ by $u$, write $\phi$ as 
$\phi(u,\bar{v}\bar{y})$, where the parameter variables $\bar{v}\bar{y}$ consist of $\bar{v}$ from the sort $\mathcal{V}$ and $\bar{y}$ from the sort $\mathcal{K}$. By Lemma~\ref{kuz-pierce}, we can write $\phi$ as a disjunction of contradictory conjunctions of (possibly negated) $\theta$-formulas, quantifier-free formulas, and field formulas. Since we can sum the cardinalities of disjoint definable sets and add the corresponding polynomials, and since the definability condition lifts to the disjunction, we may suppose that $\phi$ itself is such a conjunction. 

Thus we assume that $\phi$ is $\displaystyle{\bigwedge_{i=1}^n \psi_i}$ where each $\psi_i$ is a possibly negated $\theta$-formula, a possibly negated formula of form $t(u,\bar{v},\bar{w})=0$, or a field formula. Since $u$ does not occur in a field formula, we may ignore these. Likewise, a term equality is satisfied by one or all elements of $V$, so we may ignore these as well. Hence, we may suppose that each $\psi_i$ is a $\theta$-formula or its negation. 

We focus on one such formula, say 
$\displaystyle{\theta_n\left(z_1u+\sum_{i=1}^r z_{1j}v_j,\ldots, z_nu+\sum_{i=1}^nz_{nj}v_j\right)}$, where the $z_i$ and  $z_{ij}$ are $L_R$-terms, that is, polynomials in the field variables $\bar{y}$. Dividing out by non-zero $z_i$, which we may do by increasing the initial set of disjuncts, and collecting terms, we can write this formula in the form $\theta_n(u+w_1,\ldots,u+w_m,w_1',\ldots,w'_{m'})$, where $m+m'=n$ and the $w_i$ and $w_i'$ are terms in $\bar{v}$ and $\bar{y}$. In what follows, the notation $\langle v_1,\ldots , v_k\rangle$ for vectors $v_1,\ldots ,v_k$ denotes the $F\/$-span of $v_1,\ldots ,v_k$. We then have 
\begin{eqnarray*}
&&M\models \theta_n(u+w_1,\ldots,u+w_m,w_1',\ldots,w'_{m'})\\
&\Leftrightarrow& (\forall \bar{a}\in F^m) (\forall \bar{b}\in F^{m'})\left[\left(\sum_{i=1}^ma_i(u+w_i)+\sum_{i=1}^{m'}b_iw_i'=0\right)\to \left(\bigwedge_{i=1}^m a_i=0\wedge \bigwedge_{i=1}^{m'}b_i=0\right)\right]\\
&\Leftrightarrow& (u\not\in \langle \bar{w},\bar{w}'\rangle_F) \wedge \left(\mbox{there is no non-trivial relation~} \sum_{i=1}^m c_iw_i+\sum_{i=1}^{m'}d_iw_i'=0 \mbox{~with~} \sum_{i=1}^m c_i=0\right)\\
&&\mbox{{\bf or} }  \left( \mbox{$\bar w, \bar w'$ are linearly independent } \wedge u=\sum_{i=1}^m c_iw_i+\sum_{i=1}^{m'}d_iw_i'\mbox{ with } \sum_{i=1}^m c_i\not=-1\right) .
\end{eqnarray*}
The first disjunct in the last equivalence yields a definable set that has cardinality 
$|V|-|F|^{\dim\langle\bar{w}\bar{w}'\rangle}$ 
or 0. The second disjunct, which follows by linear algebra (we omit the details), gives a definable set of size 
$|F|^{m+m'}-|F|^{m+m' -1}$. In both cases, the conditions on $\bar{w}\bar{w}'$ determine the size of the definable set, which in turn are determined entirely by the {\em parameters} $\bar{v}\bar{y}$.

By breaking up the formula $\phi$ according to conditions on the parameters, we reduce to the case where $\phi$ is a conjunction of $\theta$-formulas (possibly negated) in terms in $u$, $\bar{v}$, and $\bar{y}$, 
expressing that $u\in (U_1\cap \ldots \cap U_l)\setminus (V_1\cup\ldots \cup V_k)$ where the $U_i$ and $V_i$ are cosets of subspaces spanned by terms in the parameters, possibly with $l=0$. Such a set has cardinality given by one of finitely many polynomials in $|V|, |F|$, of the form $|V|+p(|F|)$, or $p(|F|)$. The parameters corresponding to each polynomial are uniformly definable in the parameters of the formula.

\smallskip

\noindent {\em Case 2.}\enspace We now suppose that $x$ lies in the field sort. Again, by breaking $\phi$ into a disjunction of conjunctions, we may assume that $\phi$ is a conjunction of field formulas, along with possibly negated equations and $\theta\/$-formulas in the vector space sort. By breaking $\phi$ into more disjunctions, we shall reduce to the case where each conjunct is a field formula. In this case, the fact that finite fields form a 1-dimensional asymptotic class (the content of the main theorem of \cite{cdm}) is applicable, with the polynomials being monomials of the form $\mu$ or $\mu|F|$, where $\mu\in {\mathbb Q}$. The definability clause (concerning the formulas $\phi_p$) is also easily checked. 

Equations in the vector space sort take the form
$\displaystyle{\sum_{i=1}^t p_i(x,\bar{y})v_i=0}$, where each $p_i$ is a polynomial. If the $v_i$ are linearly independent---a condition on the parameters---this is equivalent to $\displaystyle{\bigwedge_{i=1}^t p_i(x,\bar{y})=0}$, a field formula. If $v_1,\ldots,v_{t'}$ are linearly independent, and $v_{t'+1},\ldots,v_t$ are in the span of $v_1,\ldots,v_t$, then the corresponding scalars exhibiting this are definable in $\bar{v}$, and the original equation becomes equivalent to a field formula in $x,\bar{y}$, and these scalars.

Lastly, we consider $\theta$-formulas, which (possibly negated) have form
$$\theta_n\left(\sum_{j=1}^t p_{1j}(x,\bar{y})v_j,\ldots,\sum_{j=1}^tp_{nj}(x,\bar{y})v_j\right).$$
If $v_1,\ldots,v_t$ are linearly independent, this formula is equivalent to a field condition on $x\bar{y}$ (that the matrix 
$(p_{ij}(x,\bar{y}))$ has rank $n$). And again, if
$v_1,\ldots,v_{t'}$ are linearly independent, and $v_{t'+1},\ldots,v_t$ lie in the span of $v_1,\ldots,v_{t'}$, then $\theta$ is equivalent to a field condition in $x,\bar{y}$, and the scalars exhibiting linear dependence, which, as before, are definable in the parameters $\bar{v}$. Lemma~\ref{r=1} is now proved. 
\edem

{\em Proof of Theorem~\ref{cdmvs}.} We prove Theorem~\ref{cdmvs}  by induction on $r$, mimicking the proof of \cite[Lemma 2.2]{elwes}, with many details omitted. By Lemma~\ref{r=1} the induction starts. Assume the proposition holds for $r$, and consider a formula $\phi(\bar{x},\bar{y})$ where $l(\bar{x})=r+1$. Put
$\bar{x}=z\bar{x}'$. By the case $r=1$, there is a finite set $E$ of polynomials $p({\bf V},{\bf F})\in {\mathbb Q}[{\bf V},{\bf F}]$ such that 
Theorem~\ref{cdmvs} holds for the formula $\phi(z,\bar{x}'\bar{y})$, and there are corresponding formulas $\phi_p$ for each $p\in E$. Put
 $E=\{p_1,\ldots,p_t\}$, and $\phi_i=\phi_{p_i}$ for each $i=1,\ldots,t$. 
 
 By the induction assumption, for each $p_i\in E$, there is a finite set $E_i$ of 
polynomials $q({\bf V},{\bf F})\in {\mathbb Q}[{\bf V},{\bf F}]$ such that for all $M=(V,F)\in \C$ and $\bar{a}\in M^s$, there is $q\in E_i$ such that
\begin{equation*}
\big||\phi_i(M^r,\bar{a})| - q(|V|,|F|)\big|= o(q(|V|,|F|)), 
\end{equation*}
along with a corresponding formula $(\phi_i)_{q}(\bar{y})$. Put $E_i=\{q_{ij}
:j=1,\ldots,r_i\}$, and $\phi_{ij}:=(\phi_i)_{q_{ij}}$ for each $i,j$.  Observe for each $M\in \C$ 
and $\bar{a}\in M^s$ that there is a unique function
$h:\{1,\ldots,t\}\to \omega$ such that for each $i=1,\ldots,t$, we have $M\models \phi_{ih(i)}(\bar{a})$. Also, for each $i$, we have $h(i)\in \{1,\ldots,r_i\}$, so the set of all such $h$ is finite. 

Now fix $M\in \C$ and $\bar{a}\in M^s$, and let $h:\{1,\ldots,t\}\to \omega$ be the corresponding function as above, so that 
$M\models \phi_{ih(i)}(\bar{a})$ for each $i=1,\ldots,t$. Then for each $i=1,\ldots,t$ we have 
$$\big|\phi_i(M^r,\bar{a})-q_{ih(i)}(|V|,|F|)\big|=o(q_{ih(i)}(|V|,|F|)),$$
and for each $\bar{b}\in \phi_i(M^r,\bar{a})$, we have
$$\big|\phi(M,\bar{b},\bar{a})-p_i(|V|,|F|)\big|=o(p_i(|V|,|F|)).$$
Let $\displaystyle{P({\bf V},{\bf F})=\sum_{i=1}^t p_i({\bf V},{\bf F}).q_{ih(i)}({\bf V},{\bf F})}.$
 Then 
$$\big|\phi(M^{r+1},\bar{a})-P(|V|,|F|)\big|=o(P(|V|,|F|)).$$
This completes the proof of Theorem~\ref{cdmvs} $\Box$

\begin{remark} \rm
The proof of Theorem~\ref{cdmvs} provides somewhat more information than stated. First, for a formula $\phi(u,\bar{y})$ where 
$u$ ranges 
through the vector space sort, there is a finite set $E$ of pairs $(k,p({\bf F}))$ where $k\in \{0,1\}$, $\mu \in {\mathbb Q}$,  and $p({\bf F})\in {\mathbb Q}[{\bf F}]$ such that for 
all $M=(V,F)\in \C$ and $\bar{a}\in M^s$, there is $(k,p({\bf F}))\in E$ such that $\phi(M,\bar{a})$ has size {\em exactly} $k|V|+p(|F|)$. Likewise, 
for any such $\phi(x,\bar{y})$ where $x$ ranges through the field sort, there is a constant $C$ and a finite set $E\subset {\mathbb Q}^{>0}$ 
such that for any $M=(V,F)\in \C$ and $\bar{a}\in M^s$, either $|\phi(M,\bar{a})|\leq C$ or there is $\mu \in E$ such 
that $\big||\phi(M,\bar{a})|-\mu|F|\big|\leq C|F|^{\frac{1}{2}}$. 

More generally, for any formula $\phi(\bar{x},\bar{y})$, there is a finite set $E$ of polynomials   
$P({\bf V},{\bf F})\in {\mathbb Q}[{\bf V},{\bf F}]$, each of the form $\prod_{i=1}^d (k_i{\bf V}+p_i({\bf F})).\mu{\bf F}^e$, with 
$k_i\in \{0,1\}$, $d,e< \omega$ and $p_i({\bf F})\in {\mathbb Q}[{\bf F}]$, and a constant $C>0$, such that for all $M=(V,F)\in \C$ and $\bar{a}\in M^s$
there is $P({\bf V},{\bf F})$ as above such that 
$$\big||\phi(M^r,\bar{a})|-P(|V|,|F|)\big|\leq C\prod_{i=1}^d (k_i|V|+p_i(|F|)).|F|^{e-\frac{1}{2}}.$$
The corresponding definability clauses also hold. 
We omit the details---compare the proof of \cite[Theorem 2.1]{ms}. 
\end{remark}

Now, let $T_{pvf}$ be the theory of all finite-dimensional vector spaces over finite fields, and let $K$ be a pseudofinite field.
Then $T_{vs}(K)$ is a completion of $T_{pvf}$. Recall that a complete theory $T$ is said to be {\em near model complete} if, modulo the theory, every formula is equivalent to a Boolean combination of existential formulas in the same variables.

\begin{proposition} \label{nmc}
(i) The theory $T_{vs}(K)$ is near model complete.

(ii) $T_{vs}(K)$ is a supersimple theory such that the vector space sort has rank $\omega$ and the field sort has rank 1. If $(V_0,K_0)\models T_{vs}(K)$ then  $K_0$ is stably embedded, with the $\emptyset$-definable subsets of $K_0^n$ given just by the structure of $K_0^n$ in the language of rings.  
\end{proposition}

\bdem
 Assertion~(i) follows immediately from Lemma~\ref{kuz-pierce} together with the near model-completeness of any complete theory of pseudofinite fields;
the latter follows from Kiefe \cite{kiefe} -- see also Chatzidakis \cite{chatzidakis}. In fact, the near model completeness assertion can be strengthened, and is uniform across all theories of pseudofinite fields.

For (ii), observe that $T_{vs}(F)$ is interpretable in the SU-rank $\omega$ theory ACFA. For if $(K,\sigma)\models$ ACFA, then $\Fix(\sigma)$ is a rank 1 pseudofinite field, and $K$ is an infinite degree extension of $\Fix(\sigma)$ so may be viewed as an infinite dimensional vector space over $\Fix(\sigma)$. Furthermore,  any theory of pseudofinite fields occurs as the fixed field theory in some completion of ACFA; for example, if $F$ if the ultraproduct
$\prod_{i\in\omega}\Ff_q/{\mathcal{U}}$, then by the main theorem of \cite{hrushfrob}, the difference field $\prod_{i\in\omega}(\Ff_q^{{\rm alg}},x\mapsto x^q)/{\mathcal{U}}$ is a model of ACFA and has fixed field $F$. The assertions about the induced structure on $K_0$ can be derived directly from (i), or from the corresponding statements about ${\rm Fix}(\sigma)$ in $(K,\sigma)\models$ ACFA -- see e.g. \cite[Proposition 5.3]{chatznd}.
\edem

For interest, we make some further observations on $T_{vs}$.
Recall \cite[Defnition 2.7]{hp} that a global $A$-invariant type $p(\bar{x})$ over a large saturated model $\Uu$ is said to have {\em NIP} if every Morley sequence $(\bar{b}_i:i <\omega)$ in $p$ over $A$---which has a uniquely determined complete type over 
$A$--- has the property that for every formula $\phi(\bar{x},\bar{y})$ there is $n_\phi<\omega$  such that for any $\bar{c}$, there are at most $n_\phi$ alternations of truth values of $\phi(\bar{b}_i,\bar{c})$ as $i$ increases.
Also, following \cite{pillay1}, a global type $p\in S(\Uu)$ is {\em generically stable\/} if it is invariant over some small set 
$A \subset \Uu$ and if for some (every) Morley sequence $(\bar{a}_i:i<\lambda)$ in $p$ over $A$ 
 and every formula $\phi(\bar{x})$ (not necessarily over $A$)
$\{i\in \lambda: \phi(\bar{a}_i)\}$ is finite or cofinite in $\lambda$. This definition differs slightly from that in \cite{hp}, but agrees for types with NIP.

\begin{proposition} \label{genstab}
There is a unique complete 1-type $p$ over $V(\Uu)$ containing all formulas of the form
$v\not\in \la u_1,\ldots,u_r\ra$, where $u_1,\ldots,u_r\in \Uu$. The type $p$ is invariant over $\emptyset$, is NIP in the sense of \cite[Remark 2.7]{hp}, and is generically stable.

\end{proposition}

\bdem
The uniqueness and invariance are clear. For uniqueness, suppose that $u,u'$ both satisfy the prescribed formulas over $\Uu$, and let $\Uu'$ be an elementary extension of $\Uu$ containing $u,u'$. Then there is 
$g\in \Aut(\Uu')$ fixing pointwise $\Uu$ and $K(\Uu')$ (so $g$ is linear over $K(\Uu')$) with $g(u)=u'$. 

To see that $p$ is NIP, let $(v_i:i\in \omega)$ be a Morley sequence over $\emptyset$, that is, a sequence of linearly independent vectors, and let $\phi(x,\bar{y}\bar{z})$ be a formula.
Let $\bar{c}, \bar{a}$ be tuples from $\Uu(V)$ and $\Uu(K)$, respectively. At most $l(\bar{c})$ of the $v_i$ are in $\la \bar{c}\ra$, and it follows that $\phi(v_i,\bar{c},\bar{a})$ has at most $l(\bar{c})+1$ alternations of truth value.

To see that $p$ is generically stable, by \cite[Proposition 3.1, Remark 3.3(iii)]{hp}, we must show that any Morley sequence in $p$ over any parameter set $A$ is totally indiscernible. This is immediate. 
\edem

\subsection{Ultraproducts of finite homocyclic groups}\label{homocyclic}
We give here an exact---rather than just asymptotic---uniformity result on the cardinalities of definable sets in homocyclic $p$-groups. Here, a {\em homocyclic group\/} is a direct sum of isomorphic cyclic $p$-groups. In this subsection, for a prime $p$ we denote by $\C_p$ the set of all finite groups $(\Z/p^n\Z)^m$ (for
 $m,n \in {\mathbb N}$). Put $\C:=\bigcup(\C_p: p \mbox{~prime})$.

\begin{theorem}\label{abelian}
Let $p$ be prime. Then any infinite ultraproduct $G$ of groups  $((\Z/p^n\Z)^n : n<\omega)$   is stable but not superstable, and  satisfies ($A$) and
 ($\DC_L$). 
 \end{theorem}
 
 We first prove the following more general result, Proposition~\ref{abelianexact}, analogous to Theorem~\ref{cdmvs}. By  the classical elimination theory for abelian groups of Szmielew (see \cite[Theorem A.2.2]{hodges}), modulo the theory of abelian groups, every formula $\phi(\bar{x},\bar{y})$ is equivalent to a boolean combination of formulas of form $t(\bar{x},\bar{y})=0$ and $p^\ell|t(\bar{x},\bar{y})$, where $t$ is a term in the language of groups, and $p$ is a prime; we shall say that such a formula is in {\em standard form}. For an abelian group $G$, we write  $G[p^k]$ for the subgroup $\{g\in G: p^kg=0\}$. For nonnegative integers $d,k$, let $S(d,k)$ be the set of functions of the form $P(X,u,v)=\sum_{i=0}^k\sum_{j=-kd}^{kd} c_{ij}X^{u(iv+j)}$, where 
 $c_{ij} \in {\mathbb Z}$ for all $0\leq i\leq k$ and $-kd\leq j\leq d$.
 
 \begin{proposition}\label{abelianexact}
Let $\phi(\bar{x},\bar{y})$ be a formula in the language of groups in standard form. Let $d$ be the greatest integer $l$ such that for some prime $p$, either some subformula $p^l|t(\bar{x},\bar{y})$ occurs in $\phi$ or some term $t(\bar{x},\bar{y})$ occurring in 
$\phi$ has a coefficient divisible by $p^l$. Then 
\begin{itemize} 
 \item[(i)] There  is a finite subset $F=F(\phi)$ of $S(d,r)$ (where $r=l(\bar{x})$) such that for each $G=({\mathbb Z}/p^n{\mathbb Z})^m \in \C$ and $\bar{a}\in G^r$, there is $P(X,u,v)=\sum_{i=0}^k\sum_{j=-kd}^{kd} c_{ij}X^{u(iv+j)} \in F$ with $c_{ij}=0$ whenever $in+j<0$, such that $|\phi(G^r,\bar{a})|=P(p,m,n)$.
 
 \item[(ii)] For each such function $P\in F$ there is a formula $\phi_P$ such that for each $G=(\Z/p^n\Z)^m\in \C$ and $\bar{a}\in G^s$, we have
  $G\models \phi_P(\bar{a})$ if and only if $|\phi(G^r,\bar{a})|=P(p,m,n)$.
\end{itemize} 
  
In addition, the same value of $d$ suffices if $\phi(\bar{x},\bar{y})$ is replaced by a formula of form 
$\bigwedge_{i=1}^s\phi(\bar{x},\bar{y}_i)$, for a fixed $s$.
 
 \end{proposition}
 
\bdem The proof is by induction on $r=l(\bar{x})$, using a fibering argument like that in the proof of Theorem~\ref{cdmvs}. To start the induction consider a formula $\phi(x,\bar{y})$ in standard form in the group $G=(\Z/p^n\Z)^m$. This formula is a boolean combination of formulas of form $t(x,\bar{y})=0$ and $q^\ell|t(x,\bar{y})$, where $t$ is a term in the language of groups, and $q$ is a prime. Clearly, we may assume $q=p$, since if $(q,p)=1$ then every element of $G$ is $q^l$-divisible. Also, as in the proof of Lemma~\ref{r=1}, we may assume $\phi$ is a conjunction of such formulas and their negations. A formula of form $t(x,\bar{a})=0$ either defines $\emptyset$ or $G$, or  has 
the same solution set as a formula of the form $p^lx=a'$ (where $l\leq d$), and so defines in $G$ a coset of the subgroup $G[p^l]$ of order $p^{lm}$.  

Consider now a formula $p^\ell|t(x,\bar{a})$. First observe that the formula $p^\ell|x$ has exactly $(p^{n-\ell})^m$ solutions in $G$. Now suppose $t(x,\bar{y})$ has form 
$\displaystyle{kx+\sum_{i=1}^s n_iy_i}$. Let $k=p^j\cdot k'$ where $(p,k')=1$. As the map $z\mapsto p^jz$ has kernel of size 
$(p^j)^m$, the formula
$p^\ell|t(x,\bar{a})$ has solution set of size  $p^{jm}\cdot (p^{n-\ell})^m=p^{m(n+j-\ell)}$ with $\ell\geq j$, or $p^{nm}$ solutions, or no solutions. In particular, it defines $\emptyset$ or a coset of the subgroup $p^iG$ of $G$ of order $p^{(n-i)m}$ for some $i\leq d$. Thus, it has exactly $p^{(n-i)m}$ solutions, with $n\geq i$, or no solutions, where $i$ is determined just by the original formula $p^\ell|t(x,\bar{y})$.

To complete the proof for $r=1$, consider an arbitrary conjunction of such formulas or their negations. A finite conjunction of such formulas again has solution set $\emptyset$ or a coset of a subgroup of order 
$p^{im}$ or $p^{(n-i)m}$ for some $i\leq d$ determined by the conjunctions. We use here that $G[p]<G[p^2]<\ldots G[p^u]=p^{n-u}G<\ldots pG<G$ and that each conjunction, if consistent, defines a coset of some group in this chain. It follows, using Inclusion-Exclusion,  that for a formula $\phi(x,\bar{y})$ in standard form, there is
a finite set $E$ of tuples $\bar{e}\bar{e}'$ of integers, where $\bar{e}=(e_0,\ldots,e_d)$ and $\bar{e}'=(e_0',\ldots,e_d')$
 such that for all $\bar{a}\in G^s$, there is $\bar{e}\bar{e}'\in E$ with
 $\displaystyle{|\phi(G,\bar{a})|=\sum_{j=0}^d e_jp^{jm}+ e_j'p^{(n-j)m}}$. Furthermore, 
for each $\bar{e}\bar{e}'\in E$ there is a formula $\phi_{\bar{e}\bar{e}'}$ defining the corresponding set of $\bar{a}\in G^s$, uniformly as $G$ ranges through $\C'$. Putting $c_{0j}=e_j$ for $0\leq j\leq d$, $c_{0j}=0$ for $-d\leq j<0$,  and $c_{1j}=e_{-j}'$ for $j$ with $-d<j\leq 0$ and $c_{1j}=0$ for $0<j\leq d$, we see that 
$|\phi(G,\bar{a})|=\sum_{i=0}^1\sum_{j=-d}^d c_{ij}p^{m(in+j)}$, as required for $r=1$.

The proposition for formulas $\phi(\bar{x},\bar{y})$ now follows by a standard fibering argument as in the proof of 
Theorem~\ref{cdmvs}. As the notation is intricate, we provide details.

Let $\phi(\bar{x},\bar{y})$ be a formula with $l(\bar{x})=r+1$, and put
$\bar{x}=z\bar{x}'$. By the case $r=1$, there are $t<\omega$ and a finite subset $D=\{P_1,\ldots,P_t\}$ of $S(d,1)$, such that for any $\bar{b}\bar{a}\in M^{r+s}$ there is $i\in \{1,\ldots,t\}$ such that
$|\phi(G,\bar{b},\bar{a})|=P_i(p,m,n)$. Furthermore, there is for each $i=1,\ldots,t$ a further formula $\phi_i(\bar{x}',\bar{y})$ defining the set of such $\bar{b}\bar{a}$. 

By induction, for each $i\in \{1,\ldots,t\}$ there is a finite set $D_i$ of $S(d,r)$ such that for any $G\in \C$ and $\bar{a}\in G^s$, there is $Q(X,u,v) \in D_i$ such that $|\phi_i(G^r,\bar{a})|=Q(p,m,n)$. Let $D_i=\{Q_{i1},\ldots,Q_{i,r_i}\}$ for each $i$, and let $\phi_{ij}(\bar{y})$ be the corresponding formula defining the set of such $\bar{a}$. As with Theorem~\ref{cdmvs}, for each $G\in \C$ and $\bar{a}\in G^r$, there is a unique function $h:\{1,\ldots,t\}\to \omega$ such that $G\models \phi_{ih(i)}(\bar{a})$ for each $i$. We then have, for each $i=1,\ldots,t$,  
$$|\phi_i(G^r,\bar{a})| =Q_{ih(i)}(p,m,n),$$
and, for each $\bar{b}\in \phi_i(G^r,\bar{a})$,
$$|\phi(G,\bar{b},\bar{a})|=P_i(p,m,n).$$
Put $R(X,u,v)=\sum_{i=1}^t P_i(X,u,v) Q_{ih(i)}(X,u,v)$. Then $|\phi(G^{r+1},\bar{a})|=R(p,m,n)$, and $R\in S(d,r+1)$, that is, $R$ has the required form. As there are finitely many such functions $h$, the set $F$ of all possible functions $R$ is also finite.

The final assertion follows immediately from the way that $d$ is defined.
\edem

\medskip

\noindent{\em Proof of Theorem~\ref{abelian}.}\enspace The fact that the ultraproduct is stable but not superstable follows immediately from the fact that it is an abelian group with a descending chain of definable subgroups (of the form $G> pG> p^2G)>\ldots$) each of infinite index in it predecessor. 

The condition ($\DC_L$) follows easily from Proposition~\ref{abelianexact}(ii); indeed, two definable sets $X$ and $Y$ satisfy $\delta(X)=\delta(Y)$ precisely if, on a set in the ultrafilter, the corresponding definable sets  have cardinalities which are polynomials of the same degree in $p$, and this is a definable condition. It follows from Theorem~\ref{supersimple2} that (SA) does not hold. Finally, condition (A) follows from the final assertion of Proposition~\ref{abelianexact}. For given $G=(\Z/p^n\Z)^n$, any positive $\phi$-formula defines a set of size
$\sum_{i=0}^r\sum_{j=-rd}^{rd} c'_{ij}p^{n(in+j)}$, where $r,d$ depend just on $\phi$ and the $c_{ij}'$ just on the number of conjuncts, not on $m$.
$\Box$

\section{Further Properties of ($A$), ($SA$), ($\DC_L$), and ($FMV$)} \label{furtherprop}

We first consider a number of technical questions around our conditions (A), (SA), $(\FMV)$, and $(\DC_L)$: the independence theorem and stable formulas; 1-variable criteria for the conditions; transferrability to $M^{\eq}$. We also explore consequences of assuming $(\FMV)$ and $(\DC_L)$, obtaining a pregeometry under an extra hypothesis.

\subsection{The Independence Theorem}

We first observe that, under the strong hypotheses ($\DC_L$) and ($SA$), the Independence Theorem has the following translation.

\begin{proposition} \label{independence}
Assume that $M$ satisfies $(\DC_L)$ and (SA), let $E=\acl^{\eq}(E)$ be a countable parameter set, and let $P_1,P_2,P_3$ be the solution sets 
in $M$ of 1-types $p_1,p_2,p_3$ over $E$. For  $i<j$ with $i,j\in\{1,2,3\}$ let $q_{ij}(x_i,x_j)$ be a 2-type over $E$ extending $p_i(x_i) \cup p_j(x_j)$, and
 let $Q_{ij}$ be the set of realizations in $M$ of $q_{ij}$. Let $\gamma_i:=\delta(P_i)$ and suppose that $\delta(Q_{ij})=\gamma_i+\gamma_j$ for each $i<j$. Then there is a 3-type $r(x_1,x_2,x_3)$ over $E$ extending $q_{12}(x_1,x_2)\cup q_{13}(x_1,x_3)\cup q_{23}(x_2,x_3)$ with $\delta(r)=\gamma_1+\gamma_2+\gamma_3$.
 \end{proposition}
 
\bdem By Theorem~\ref{low2}(i) and \cite[Theorem 6.4.6]{wagner}, strong type and Lascar strong type coincide in $T=\Th(M)$. The result thus follows directly from the usual Independence Theorem for Lascar strong types (see \cite[Theorem 5.8]{kp2}, or \cite[Theorem 2.5.20]{wagner}), via Theorem~\ref{forking2} and Lemma~\ref{fiber}. 
\edem
 
 We would like to prove a version of the Independence Theorem in the manner of the proof of Proposition 8.4.3 of \cite{ch}, that is, based directly on counting arguments rather than quoting results for simple theories. In particular, we ask:
 
\begin{question} \rm
Does a version of Proposition~\ref{independence} hold just under the assumptions $(\DC_L)$ and $(\FMV)$?
\end{question}
 
In this direction, we make two observations---Lemmas~\ref{udistable} and ~\ref{udistable2}---both of which are standard.
 
 \begin{lemma}\label{udistable}
Assume that $M$ satisfies ($\DC_L$), let $Z\subset M^t$ be definable with $\delta(Z)=\delta_0$, and let $\phi(\bar{x},\bar{z})$ and $\psi(\bar{y},\bar{z})$ be formulas implying $\bar{z}\in Z$, with $l(\bar{x})=r$ and $l(\bar{y})=s$. Let $\theta(\bar{x},\bar{y})$ be the formula, given by ($\DC_L$), which expresses that
$\delta(\phi(\bar{x},Z)\wedge \psi(\bar{y},Z))<\delta_0$. Then $\theta$ is stable.
\end{lemma}

\bdem We essentially repeat the argument of  Lemma 8.4.2  of \cite{ch}. Suppose $(\bar{a}_i,\bar{b}_i)_{i<\omega}$ is an $L^+$-indiscernible sequence from $M^{r+s}$, with
$\theta(\bar{a}_i,\bar{b}_j)$ holding whenever $i<j$. It suffices to show that $\theta(\bar{a}_i,\bar{b}_i)$ holds for each $i$. Let
$Z_i:=\{\bar{z}:\phi(\bar{a}_i,\bar{z})\wedge \psi(\bar{b}_i,\bar{z})\}$. For $i<j$, as $Z_i\cap Z_j\subseteq \phi(\bar{a}_i,Z)\cap \psi(\bar{b}_j,Z)$, 
we have $\delta(Z_i\cap Z_j)<\delta_0$. It follows by Remark~\ref{kintersections2} (or can be derived from Proposition~\ref{kintersections}) that $\delta(Z_i)<\delta_0$ for each $i$, as required. 
\edem

In part for the next section, we recall from Section 2.1  the natural notion of measure on definable subsets of a given definable set in $M$. Given $L$-formulas $\phi(\bar{x},\bar{y})$ and $\psi(\bar{x},\bar{z})$ with $l(\bar{x})=r$, $l(\bar{y})=s$ and $l(\bar{z})=t$,  and $\bar{a}\in M^s$ and $\bar{b}\in  M^t$, we  define the normalized measure $\mu_{\psi(\ov{x},\ov{b})}(\phi(\ov{x},\ov{a}))$ with respect to $\psi(\ov{x},\ov{b})$. This gives a finitely additive real-valued probability measure on definable subsets of a given definable set. 

\bd \label{measuredef}
We say that $M$ satisfies ($\MD_L$) if for each pair of formulas $\phi(\ov{x},\ov{y})$ and $\psi(\ov{x},\ov{z})$, there is a formula
$\chi^\psi_\phi(\bar{y}_1,\bar{y}_2,\bar{z})$, with $l(\bar{y}_1)=l(\bar{y}_2)=s$ and $l(\bar{z})=t$, such that for all 
$\bar{a}_1,\bar{a}_2\in M^s$ and $\bar{b}\in M^t$, we have
$$\chi^\psi_\phi(\phi(\bar{a}_1,\bar{a}_2,\bar{b}))\Leftrightarrow [\mu_{\psi(\ov{x},\ov{b})}(\phi(\bar{x},\bar{a}_1))\leq \mu_{\psi(\ov{x},\ov{b})}(\phi(\bar{x},\bar{a}_2))].$$
\ed

Note that this condition is an $\bigwedge$-definable condition in the language $L^+$. We do not consider the corresponding notion ($\MD_{L^+}$).
It is easily seen that if ($A$) and ($\MD_L$) both hold, then  for all formulas $\psi(\bar{x},\bar{z})$ and $\phi(\bar{x},\bar{y})$, the quotient total ordering of the above preordering is finite; that is, the collection of $\phi$-definable subsets of $\psi(\bar{x},\bar{b})$ assumes finitely many measures as $\bar{z}$ varies, cf., Remark~\ref{SOP}.
 
For future reference, we explicitly record the following measure-theoretic version of Lemma~\ref{udistable}, used in the  proof of Theorem~\ref{tao} below.  It follows immediately from \cite[Proposition 2.25]{h1}.
\begin{lemma} \label{udistable2}
Assume that $M$ satisfies ($\DC_L$) and ($MD_L$), and that $D$ is a definable set in $M^t$. Let $\phi(\bar{x},\bar{z})$ and $\psi(\bar{y},\bar{z})$ be formulas which imply $\bar{z}\in D$. For some $\mu>0$, let $\theta(\bar{x},\bar{y})$ define the set of all $(\bar{a},\bar{b})$ such that $\mu_D(\phi(\bar{a},M^t)\wedge \psi(\bar{b},M^t))=\mu$. Then $\theta$ is stable. 
\end{lemma}

\subsection{1-variable criteria}

By Lemma 2.2 of \cite{elwes}, to show that a collection $\C$ of finite structures is an asymptotic class, it suffices to verify the conditions for formulas $\phi(x,\bar{y})$ where $x$ is a single variable. This is analogous to combinatorial conditions on stable and NIP formulas, and our use above of Lemma~\ref{r=1}. Below, we give a clean result for ($\DC_L$) and ($\FMV$), but have not obtained analogues for ($A$) or ($SA$).
 
\bd The conditions ($\DC_L$)($k$), ($\FMV$)($k$) are defined like ($\DC_L$), ($\FMV$) respectively, but only for 
$r=l(\bar{x})\leq k$. 
\ed

\begin{lemma}\label{weakOD}
Assume ($\FMV$)($1$) and ($\DC_L$)($1$). Then ($\DC_L$) and ($\FMV$) hold.
\end{lemma}

\bdem We show by induction on $k$ simultaneously that  ($\DC_L$)($k$) and ($\FMV$)($k$) both hold.  By our assumptions,  both assertions hold for $k=1$. 

Assume both statements hold for some $k$ and consider the formula $\phi(\bar{x},\bar{y})$ where $l(\bar{x})=k+1$. Put $\bar{x}=z\bar{x}'$. Define
$Q=\{\delta(\phi(M,\bar{b}\bar{a})):\bar{b}\bar{a}\in M^{k+s}\}$; this set is finite by ($\FMV$)($1$), so we may put 
$Q=\{\gamma_1,\ldots,\gamma_t\}$.  For each $\gamma\in Q$, let $\psi_\gamma(\bar{x}',\bar{y})$ hold if and only if
$\delta(\phi(M,\bar{x}',\bar{y}))=\gamma\,$---the formula $\psi_\gamma$ exists by ($\DC_L$)($1$) and the finiteness of $Q$. Then put
$$Q(\gamma):=\{\delta(\psi_\gamma(M^k,\bar{a})): \bar{a}\in M^s\}.$$
By the induction assumption ($\FMV$)($k$), this set is finite. 
Then for all $\bar{a}\in M^s$ we have 
\begin{equation*}
\delta(\phi(M^{k+1},\bar{a}))=\sum_{i=1}^t \gamma_i+\delta(\psi_{\gamma_i}(M^k,\bar{a})),\tag{**}
\end{equation*}
and the set of all such values is finite as each $Q(\gamma_i)$ is finite. Thus 
 ($\FMV$($k+1$) holds.
 Using the induction assumption ($\DC_L$)($k$) and ($**$), it is easy to see  ($\DC_L$)($k+1$) also holds.
 \edem

\begin{question}\rm
Are there analogues of Lemma~\ref{weakOD} for the conditions (A) and (SA), possibly local for (A)?
\end{question}

\subsection{Transferring conditions to  $M^{\eq}$}
We consider here the extent to which the conditions on which this paper has focused, (A), (SA), and $(\DC_L)$, extend to $M^{\eq}$. Analogously, it is shown in \cite{elwes} that, essentially, if $\C$ is an asymptotic class and $\C'$ is obtained from $\C$ by adding finitely many sorts from $M^{\eq}$, then $\C'$ is an asymptotic class.

\begin{proposition} \label{eq}
{\rm (i)} If $M$  has ($\DC_L$) and ($\FMV$), then $M^{\eq}$ satisfies ($\DC_L$) and ($\FMV$).
\begin{itemize}
\item[(ii)] If $M$ has ($\DC_L$) and $\Th(M)$ does not have the strict order property, then $M^{\eq}$ satisfies ($\DC_L$) and 
($\FMV$).

\item[(iii)] If $M$ satisfies ($A$) and ($\DC_L$), then $M^{\eq}$ satisfies ($A$) and ($\DC_L$).

\item[(iv)] If $M$ satisfies ($SA$) and ($\DC_L$), then $M^{eq}$ satisfies ($SA$) and ($\DC_L$).
\end{itemize}
\end{proposition}

\bdem
(i) Let $E$ and $F$ be $\emptyset$-definable equivalence relations on $M^n$ and $M^m$ respectively, and let $\phi(x,\bar{y})$ and $\psi(u,\bar{v})$ be $L^{\eq}\/$-formulas such that $x$ ranges through $M^n/E$ and $u$ through $M^m/F$. Using Lemma~\ref{weakOD}, it suffices to show that the relation $\delta(\phi(x,\bar{a}))\leq\delta(\psi(u,\bar{b}))$ is defined by some formula $\chi(\bar{a},\bar{b})$.  By ($\FMV$), the $E$ and $F$-classes take just finitely many $\delta$-values. Hence, the set of  $E$-classes or $F$-classes taking any  given $\delta$-value is $\emptyset$-definable. The result now follows easily, using 
Lemmas~\ref{WOD} and~\ref{WODapp}(i).

(ii) This is immediate from Lemma~\ref{SOP}, in conjunction with~(i). 

(iii) First observe that $M^{\eq}$ satisfies ($\DC_L$), by~(i) and Theorem~\ref{low2}. To see that $M^{\eq}$ satisfies ($A$), let $E$ be an $\emptyset$-definable equivalence relation on $M^n$, and suppose that there is a sequence of subsets $(X_i:i\in\omega)$ of $M^n/E$ such that $X_i\supset X_{i+1}$ and each $X_i$ is a conjunction of $\phi$-instances for some $L^{\eq}\/$-formula $\phi$. (For simplicity, we are handling a special case where $X_i\subset M^n/E$ for each $i$; the general argument hen $X_i\subseteq (M^n/E)^r$ is similar.) For each $i<\omega$ let $X_i'$ be the union of the $E$-classes lying in $X_i$. 
 Then $X_i'\supset X_{i+1}'$ for all $i$, and it is easily checked that there is an $L$-formula $\phi'$ such that each $X_i'$ is a conjunction of $\phi'$-instances. By ($A$), there is $t$ such that $\delta(X_i')=\delta(X_t')$ for all $i\geq t$. Also, the set of $E$-classes is uniformly definable, so by ($A$) and Lemma~\ref{SOP}, these take finitely many $\delta$-values.
 It follows, again using Lemma~\ref{WODapp}(i), that the sequence $(\delta(X_i):i\in \omega)$ takes finitely many values, as required.
 
 (iv) This is proved essentially as in (iii).
 \edem
 
 \begin{remark} \rm
It would be helpful to clarify what hypotheses are needed in Proposition~\ref{eq}. 
 In~(i) above, it seems we require some assumption in addition to ($\DC_L$) for $M$. 
 Likewise, in~(iii), we probably cannot deduce that $M^{\eq}$ satisfies ($A$) just from the assumption that $M$ satisfies (A). In these cases we have not constructed counterexamples.
 
 In (iv), ($\DC_L$) is required, that is, we cannot lift ($SA$) on its own from $M$ to $M^{\eq}$. Consider a language $L$ with a binary relation $E$ and unary relations $\{P_i:i\in \omega\}$. We can choose an increasing sequence $\delta_0<\delta_1<\ldots$ and build a family of finite structures with ultraproduct $M$ such that: $E$ is an equivalence relation on $M$; each $P_i$ is a union of $E$-classes with $P_0\supset P_1\supset\ldots$; the structure $M$ satisfies ($SA$); and, $\delta(P_i)$ takes a fixed value 
 $\epsilon$ for all $i$,  but $\delta(P_i/E)>\delta(P_{i+1}/E)$ for all $i<\omega$. For this, we arrange that the $E$-classes in $P_i\setminus P_{i+1}$ all have $\delta$-value $\delta_i$.  We omit the details. 
 
 \end{remark}

\subsection{Consequences of $(\FMV)$ and (DC$_L$)}
We assume that both ($\FMV$) and ($\DC_L$) hold throughout this subsection. Note that these assumptions hold for the examples considered in Section~4: asymptotic classes, the $2\/$-sorted vector space structures with theory $T_{vs}$ of 
Theorem~\ref{cdmvs}, and the ultraproducts of homocyclic groups of Theorem~\ref{abelian}. These all satisfy ($A$) and so have a simple theory, unlike the expansion of $T_{vs}$ by a symplectic bilinear form considered in Example~\ref{bilinear}, which is not simple but does satisfy ($\FMV$) and ($\DC_L$). 

Under ($\FMV$) and ($\DC_L$) we have additivity of $\delta\/$-dimension given by Lemma~\ref{fiber} and all properties of 
$\dind$ considered in Propositions~\ref{d-properties} and \ref{d-properties+} except for local character (which fails in the symplectic bilinear form example). Also, by Proposition~\ref{eq}(i), properties ($\FMV$) and ($\DC_L$) transfer to $M^{\eq}$, and by Lemma~\ref{weakOD} it suffices to verify them for formulas of form $\phi(x,\bar{y})$. 

Our main additional observation is 

\begin{proposition} \label{geom}
Assume that $M$ satisfies ($\FMV$) and ($\DC_L$), let $D$ be an interpretable set in $M^{\eq}$ over parameters $\bar{e}$, and put $\delta_0:=\delta(D)$. Suppose there is a proper subsemigroup $S$ of $\R^*/C$ with $\delta_0\not\in S$, such that $\delta(D')\in S$ for 
every  definable subset $D'$ of $D$ with $\delta(D')<\delta_0$. For $a\in D$ and $B\subset D$, define $a\in \cl(B)$ if and only if there is a 
$B\bar{e}$-definable subset $D'$ of $D$ containing $a$ with $\delta(D')<\delta_0$.
 Then $\cl$ defines a pregeometry on $D$.
\end{proposition}

\bdem For ease of notation we suppose that $D\subset M$ and that $D$ is $\emptyset$-definable. 
We must verify, for $A,B\subset D$, that:  
\begin{itemize}
\item[(i)] if $A\subseteq B$ then $A\subseteq \cl(A)\subseteq \cl(B)$;
\item[(ii)] if $a\in \cl(B)$ then $a\in \cl(F)$ for some finite $F\subseteq B$;
\item[(iii)] $\cl(\cl(A))=\cl(A)$;
\item[(iv)] for all $a_1,a_2\in D$, we have $a_1\in \cl(a_2 B)\setminus \cl(B)\Rightarrow a_2\in\cl(a_1 B)$.
\end{itemize}

Properties~(i) and~(ii) are immediate.  To prove~(iv)  we first note that we may assume $B$ to be countable, in which case
$a\in \cl(B)\setminus \cl(\emptyset)$ if and only if $a\ndind B$. Then, if $a_1\not\in \cl(B)$ and $a_2\not\in \cl(B a_1)$ we have $a_1\dind B$ and 
$a_2\dind Ba_1$. Thus $a_2\dind_B a_1$, whence by Proposition~\ref{d-properties+}(vi) we have $a_1\dind_B a_2$, giving  $a_1\dind Ba_2$ and finally $a_1\not\in \cl(B,a_2)$.  For~(iii), suppose that $d_1,\ldots,d_r\in \cl(A)$, and 
$c\in \cl(A d_1\cdots d_r)$. For each $i=1,\ldots,r$, there is an $A$-definable set $D_i\subset D$ with $\delta_i:=\delta(D_i)<\delta_0$ and $d_i\in D_i$. Also, there is an $A d_1\cdots d_r$-definable subset $D^*\subseteq D$ with $c\in D^*$ and $\delta^*:=\delta(D^*)<\delta$. Let $\psi(x,d_1,\ldots,d_r)$ be a formula over $A$ defining $D^*$ and, using ($\FMV$) and ($\DC_L$), let 
$\chi(x,y_1,\ldots,y_r)$ be the $A\/$-formula $\psi(x,y_1,\ldots,y_r) \wedge (\delta(\psi(M,y_1,\ldots,y_r))=\delta^*)$. Now put 
$$D':=\bigcup\{\chi(M,d_1',\ldots,d_r'), d_1'\in D_1,\ldots,d_r'\in D_r\}.$$ 
Applying Lemmas~\ref{WOD} and ~\ref{WODapp}  an easy counting argument shows that 
there is an $n<\omega$ such that $|D'|\leq n|D^*| |D_1|\cdots |D_r|$. It follows that
$\delta(D')\leq \delta^*+\delta_1+\ldots+\delta_r<\delta_0$, as required. 
\edem

We shall call a set  interpretable in $M$ {\em geometric} if it  satisfies the assumptions on $D$ in Proposition~\ref{geom}. 
Examples of geometric sets include ultraproducts of  one-dimensional asymptotic classes (where the  subsemigroup $S$ in Proposition~\ref{geom} is trivial), and   both the vector space sort and the field sort in $T_{vs}$. It would be interesting to investigate the pregeometry from the viewpoint of Zilber Trichotomy, and show, for example, that in the locally modular non-trivial case there is an infinite definable group. 

We conclude this section with the proposition below for geometric sets which are groups; the analogous results for  asymptotic classes and measurable structures are \cite[Theorems 3.12, 5.15]{ms}. The conclusion cannot be strengthened to `abelian-by-finite' since for an odd prime $p$ the class of finite extraspecial $p$-groups of exponent $p$ is a one-dimensional asymptotic class (\cite[Proposition 3.11]{ms}), and has finite-by-abelian but not abelian-by-finite ultraproducts. The result suggests that if a geometric set is a pure group, it should be one-based in a Zilber Trichotomy.

\begin{proposition} \label{geom2}
Assume that $M$ satisfies ($\FMV$) and ($\DC_L$), and let $G$ be an infinite group interpretable in $M$ such that the domain of $G$ is a geometric set. Then $G$ is finite-by-abelian-by-finite. 
\end{proposition}

\bdem By a theorem of Landau \cite{landau}, for every $k< \omega$ there are just finitely many finite groups with $k$ conjugacy classes. Hence, as $G$ is a pseudofinite group, $G$ has infinitely many conjugacy classes, say $\{C_i:i\in I\}$. The conjugacy classes are uniformly definable, so the set $\{\delta(C_i):i\in I\}$ is finite. For a conjugacy class $C_i$ and $a\in C_i$, we have 
$|G|=|C_i|\cdot |C_G(a)|$ (non-standard cardinality), so $\delta_0:=\delta(G)=\delta(C_i)+\delta(C_G(a))$. If a conjugacy class $C_i$ is infinite then $|G:C_G(a)|$ is infinite. In this case $\delta(C_G(a))<\delta(G)$, and thus, as $G$ is geometric, 
$\delta(C_i)=\delta_0$. Hence, by counting, $G$ has just finitely many infinite conjugacy classes. By ($\DC_L$) the finite conjugacy classes have bounded size, which yields that the set of finite conjugacy classes of $G$ is definable and hence its union is a definable non-trivial normal subgroup $N$ of $G$.  As $G/N$ has finitely many conjugacy classes, it is finite.

The group $N$ is a so-called 
BFC group, that is, a group whose conjugacy classes have finite bounded size. It follows by \cite[Theorem 3.1]{neumann} that its derived subgroup $N'$ is finite, that is, $G$ is finite-by-abelian-by-finite. 
\edem

\section{Applications}\label{applications}
We consider here two potential routes for applications of pseudofinite dimension. The first is to pseudofinite groups, and the second to possible generalizations of Tao's `Algebraic Regularity Lemma' \cite{tao}. Many other possible lines of application, from a different viewpoint, are described in \cite{Hru}.

\subsection{Pseudofinite dimension and groups}

We here assume that there is a group $G$ definable in $M$. Under the assumption ($A$), the entire theory of groups with simple theory is applicable. We first
give a small adaptation of some observations from  \cite[Section 4]{mt}, where applications to finite simple groups of fixed Lie rank are described. These are the subject of~\ref{wordmaps}-\ref{6.1.5}. We conclude this subsection with 
Proposition~\ref{329converse}, a partial converse for expansions of groups to Theorem~\ref{supersimple2}.

Below, if $U,V\subset G$ then $UV:=\{uv: u\in U, v\in V\}$. 

\begin{theorem} \label{wordmaps}
Assume that $M$ satisfies ($SA$) and ($\DC_L$). Suppose that $G$ is a definable group in $M$ that has no proper definable subgroups of finite index, and let $X_1,X_2,X_3$ be definable subsets of $G$ with $\delta(X_i)=\delta(G)=\delta_0$ for each $i$. Then
\begin{itemize}
\item[(i)] $X_1X_2X_3=G$, and

\item[(ii)] $\delta(G\setminus X_1X_2)<\delta_0$.
\end{itemize}
\end{theorem}

Before giving the proof, we collect some basic facts about generic types in simple theories, taken from \cite[Section 4.3]{wagner} (and originally in \cite{pillay98}). First, following \cite{wagner}, by a {\em type-definable group} $G$ we mean a type-definable set together with a definable binary operation which induces a group operation on the domain. By \cite[Theorem 5.5.4]{wagner}, if $G$ is a type-definable group in a supersimple theory, then $G$ is an intersection of definable groups. In particular, if $G$ is a type-definable subgroup of the definable group $H$, then $G$ is an intersection of definable subgroups of $H$. 

Given a type-definable group $G$ and a countable set $A$ of parameters, by $S_G(A)$ we denote the set of complete types over $A$ which contain the formula $x\in G$.  For such a group $G$ in an ambient simple theory, given a countable set $A$ of parameters, a type $p\in S_G(A)$ is called {\em left generic\/} if for all $b$ realizing a type in $S_G(A)$, and all $a\models p$ with 
$a\ind_A b$, we have $ba\ind A,b$. There is an analogous definition of `right-generic', but the two notions coincide  (\cite[Lemma 4.3.4]{wagner}), so we just call such a type {\em generic}. It is shown in \cite[p.168]{wagner} that if the ambient theory is supersimple, then $p$ is a generic type of $G$ if and only if $SU(p)=SU(G)$. Recall also that if $X$ is an $A$-definable set, then $SU(X)$ is the supremum of the ranks $SU(p)$, as $p$ ranges over types over $A$ concentrating on $X$. In a general supersimple theory, this supremum may not be realised, but if X is a group, then it is realized by any generic type (see \cite[Section 5.4]{wagner}). 


\begin{lemma} \label{SUdelta}
Assume ($SA$) and ($\DC_L$). Suppose that $G$ is a group definable in $M$ and let $X\subset G$ be definable with $\delta(X)=\delta(G)$. Then $SU(X)=SU(G)$, and $X$ realizes a generic type of $G$ (over any small parameter set). 
\end{lemma}

\bdem Adding constants if necessary, we may assume that $X$ and $G$ are $\emptyset$-definable. Let $p$ be a generic type of $G$ over 
$\emptyset$. Choose $b\in G$ realizing $p$ and $c\in Xb^{-1}$ with $c\dind b\/$; this is possible as $\delta(Xb^{-1})=\delta(X)=\delta(G)$, as per Proposition~\ref{d-properties}(i). Then $c\ind b$ by Theorem~\ref{forking2}, and  $cb\in X$. Also, $cb\ind c$ 
as $b$ is generic. Hence, as $cb$ and $b$ are interdefinable over $c$, $cb$ is generic over $c$ and  so is generic over $\emptyset$, by  \cite[Lemma 4.1.2(1), (3)]{wagner}. In particular, $SU(X)=SU(G)$.
\edem

\noindent {\em Proof of Theorem~\ref{wordmaps}}.\enspace By Theorem~\ref{supersimple2}, $\Th(M)$ is supersimple. 

(i)  We work over a countable elementary submodel $M_0\prec M$. Observe that $G_{M_0}^o$, the smallest $M_0$-type-definable subgroup of $G$ of bounded index, is equal to $G$. Indeed, as noted above, $G_{M_0}^o$ is an intersection of definable subgroups of $G$ of bounded index, and by compactness and saturation such subgroups have finite index in $G$, so equal $G$ by assumption.

By Lemma~\ref{SUdelta} and our assumption for each $i$ that $\delta(X_i)=\delta(G)$, we have $SU(X_i)=SU(G)$ for $i=1,2,3$. 
Since $G=G_{M_0}^o$, every generic type of $G$ over $M_0$ is {\em principal} (see \cite[Definition 4.4.6]{wagner}). Hence, by Proposition 4.7 (ii) of \cite{mt}---a small translation of \cite[Proposition 2.2]{psw}---if $r_1,r_2,r_3$ are generic types of $G$ over $M_0$ and $r$ is any type of $G$ over $M_0$, there are $a_i\in G$ realizing $r_i$ for $i=1,2,3$ such that $a_1a_2a_3\models r$.

For a contradiction suppose that $X_1X_2X_3\neq G$. Let $r$ be any type over $M_0$ containing the formula $x\in G\setminus X_1X_2X_3$. By Lemma~\ref{SUdelta}, For each $i$ there is a generic type $p_i$ of $G$ over $M_0$ containing the formula $x\in X_i$. By the conclusion of last paragraph, there are $a_i\models p_i$ in $G$, for each $i$, with $a_1a_2a_3\models r$. This, however, contradicts the assumption that $r$ contains the formula $x\in G\setminus X_1X_2X_3$.

(ii) The proof is a small adaptation of that of \cite[Theorem 4.8(ii)]{mt}. If the conclusion were false, then (arguing as in (i)) there would be a generic type $q$ of $G$ (over $M_0$) containing the formula  $x\in G\setminus X_1X_2$. Choosing $p_1,p_2$ as in~(i), we find by \cite[Proposition 2.2]{psw} realizations $a_1\models p_1$ and $a_2\models p_2$ with $a_1a_2\models q$, a contradiction. \hskip3.9in$\Box$

\medskip

Theorem~\ref{wordmaps}(i) has  consequences for finite groups. For example we have 

\begin{corollary} \label{npyb}
Let $\C$ be a class of finite groups (possibly with extra structure) such that all ultraproducts of $\C^+$ satisfy ($SA$) and 
($\DC_L$). Assume for each positive integer $d$ and formula $\psi(x,\bar{y})$ that there are only finitely many $G\in \C$ containing a tuple $\bar{a}$ such that $\psi(G,\bar{a})$ is a proper  subgroup of $G$ of index at most $d$. Let $N<\omega$ and let $\chi_i(x,\bar{z}_i)$ for $i=1,2,3$ be formulas. Then there is $K<\omega$ such that if $G\in \C$ with $|G|>K$ and 
$\bar{a}_i\in G^{l(\bar{z}_i)}$ with $|\chi_i(G,\bar{a}_i)|\geq \frac{1}{N}|G|$ for $i=1,2,3$, then 
$\chi_1(G,\bar{a}_1)\cdot\chi_2(G,\bar{a}_2)\cdot\chi_3(G,\bar{a}_3)=G$.
\end{corollary}

Corollary~\ref{npyb} is analogous to Corollary 1 of Nikolov-Pyber \cite{np} (see also \cite[Remark 3.3]{Hru}, which concerns Hrushovski's `coarse pseudofinite dimension'). The latter, which uses a result of Gowers \cite{gowers}, has no model-theoretic assumptions, but assumes that the groups have no non-trivial  representations of bounded finite degree. The Nikolov-Pyber theorem has the following pseudofinite consequence, noted also in \cite{Hru}. Here `internal' has the usual meaning from non-standard analysis;  an internal representation would arise as an ultraproduct of representations of the finite groups.

\begin{theorem}
Let $G$ be an infinite ultraproduct of finite groups,  with no non-trivial internal finite degree representation. Let $X_1,X_2,X_3$ be definable subsets of $G$ with $\delta(X_i)=\delta(G)$ for all $i$. Then $X_1X_2X_3=G$.
\end{theorem}

Since any family of finite simple groups of fixed Lie type is an asymptotic class, Theorem~\ref{wordmaps} (via Proposition~\ref{asymp}) has the following consequence for finite simple groups, already noted in \cite{mt} and derivable also from the Nikolov-Pyber theorem. There is a much stronger statement in \cite{lst}, where the result is proved with two words rather than three, and without the restriction on Lie type. If $w(x_1,\ldots,x_d)$ is a non-trivial word in the free group on $x_1,\ldots,x_d$ and $G$ is a group, then $w(G):=\{w(g_1,\ldots,g_d): g_1,\ldots,g_d\in G\}$.

\begin{theorem} \label{6.1.5}
Let $w_1,w_2,w_3$ be non-trivial group words. Then for any fixed Lie type $\tau$ there is $N=N(w_1,w_2,w_3,\tau)$ such that if $G$ is a finite simple group of Lie type $\tau$ and $|G|\geq N$, then $G=w_1(G)\cdot w_2(G)\cdot w_3(G)$.
\end{theorem} 

Arguments with generic types also yield the following partial converse to Theorem~\ref{supersimple2}, for expansions of groups. 

\begin{proposition} \label{329converse}
Assume that $M$ satisfies ($A$) and ($\DC_L$), has a supersimple theory,  and is an expansion of a group $G$. Then $M$ satisfies ($SA^-$).
\end{proposition}

\bdem Suppose that ($SA^-$) fails. Then there is a sequence $(B_i)_{i< \omega}$ of countable subsets of $G$, with 
$B_i\subset B_{i+1}$ for each $i$ and, for $B:=\bigcup(B_i:i\in \omega)$, a type $p\in S_G(B)$, such that if $p_i:=p|B_i$ for each $i$, then $\delta(p_i)>\delta(p_{i+1})$ for each $i$. We shall define inductively a sequence of groups 
$G=G_{-1}\geq G_0 \geq G_1 \geq\ldots$, such that for $i\geq 0$ the group $G_i$ is type-definable over $B_i$ and has $p_i$ as a generic type. To start, put $G_{-1}=G$, and inductively, assuming $G_{i}$ has been defined, put

\begin{eqnarray*}
{\rm St}(p_{i+1}):= &&\{g\in G_i: p_{i+1}\cup gp_{i+1} \mbox{~does not fork over~} B_{i+1}\}\\
= &&\{g\in G_i: \mbox{~there is~}y\models p_{i+1} \mbox{~with~} g\ind_{B_{i+1}} y \mbox{~and~} gy\models p_{i+1}\},
\end{eqnarray*}
and put $G_{i+1}:={\rm St}(p_{i+1}).{\rm St}(p_{i+1})$. 
Here we follow  \cite{pillay98} (see also \cite[Definition 4.5.1]{wagner}). 
By \cite{pillay98} or \cite[p.122]{wagner}, $G_{i+1}$ is a $B_{i+1}$-type-definable subgroup of $G_{i}$ and $p_{i+1}$ is a generic type of $G_{i+1}$. 

We claim that $\delta(p_i)=\delta(G_i)$ for each $i\geq 0$. Indeed, as in Proposition~\ref{d-properties} (i), there is 
 $g\in G_i$ with $\delta(g/B_i)=\delta(G_i)$. Pick 
 $a\models p_i$ with
$a\ind_{B_i} g$. Then 
$ga\ind_{B_i} g$, as $p_i$ is a generic type of $G_i$. Applying Theorem~\ref{forking2}(i) three times we have 
$$\delta(a/B_i)=\delta(a/g,B_i)=\delta(ga/g,B_i)=\delta(ga/B_i)\geq \delta(ga/a,B_i)=\delta(g/a,B_i)=\delta(g/B_i).$$  
Thus $\delta(a/B_i)= \delta(g/B_i)$, yielding the claim. 

It follows that $(G_i)_{i<\omega}$ is a decreasing sequence of type-definable subgroups of $G$ with $\delta(G_{i+1})<\delta(G_i)$ for each $i<\omega$. We now assert that there is a decreasing sequence $(H_i)_{i<\omega}$ of {\em definable} subgroups of $G$ with $G_i\leq H_i$ and $\delta(H_{i+1})<\delta(H_i)$ for each $i$. 
To see this, suppose that $H_0>\ldots>H_n$ have been constructed satisfying these conditions. By \cite[Theorem 5.5.4]{wagner}, $G_{n+1}$ is the intersection of a family $(K_i:i\in I)$ of {\em definable} subgroups of $G$. As $\delta(\bigcap(K_i:i\in I))<\delta(G_n)$,  there is a definable set $D \supseteq \bigcap(K_i:i\in I)$ with $\delta(D)<\delta(G_n)$. By compactness and saturation there is finite $I_0\subset I$ such that
$D\supseteq\bigcap(K_i:i\in I_0)$.  Putting $H_{n+1}:=H_n\cap \bigcap(K_i:i\in I_0)$, we find $\delta(H_{n+1})<\delta(H_n)$.

By Lemma~\ref{basics}(ii), $|H_i :H_{i+1}|$ is infinite for each $i$. For $i<\omega$, let $c_i$ be a canonical parameter in 
$M^{\eq}$ for $H_i$. Let $a$ be generic in $\bigcap(H_i: i<\omega)$. Then $a\nind_{c_0\ldots c_i} c_{i+1}$ for each $i$, contradicting supersimplicity. 
\edem

\subsection{Tao's Algebraic Regularity Lemma}

We here give a generalization of Tao's  Algebraic Regularity Lemma, proved in 
\cite{tao} with a remarkable application to expansion properties for polynomials. No new ideas are involved in our treatment---it is a routine application of methods of Pillay and Starchenko
 \cite{ps}, combined with the argument from \cite{tao} to deduce his Lemma~5 from his Proposition~27. In unpublished work, Hrushovski gives a rather stronger generalization. We omit the details. Below, `complexity' refers to the length of a formula and 
 ($\MD_L$) is as in Definition~\ref{measuredef}.

\begin{theorem} \label{tao}
Let $\C$ be a class of finite $L$-structures, and assume that every infinite ultraproduct of members of $\C$ satisfies ($SA$), 
($\DC_L$), and ($\MD_L$). Then for every $N\in {\mathbb N}^{>0}$ there is $C=C_N\in {\mathbb N}^{>0}$ such that: whenever $M\in \C$ has cardinality greater than $C$,
$V$ and $W$ are non-empty  sets in cartesian powers of  $M$,  and $E\subseteq V \times W$, with $V,W$ and $E$ all definable of complexity at most $N$, there are partitions $V=V_1\cup \ldots \cup V_a$ and $W=W_1\cup \ldots \cup W_b$ into definable sets of complexity at most $C$, with:
\begin{itemize}
\item[(1)] for all $i=1,\ldots,a$ and $j=1,\ldots,b$, we have $|V_i|\geq |V|/C$ and $|W_j|\geq |W|/C$, and 

\item[(2)] for all $i,j$, and sets $A\subset V_i$ and $B\subset W_j$, we have
$$\big||E \cap (A\times B)|-d_{ij}|A||B|\big|=o(|V_i||W_j|),$$
where $d_{ij}=|E\cap(V_i\times W_j)|/|V_i||W_j|$.
\end{itemize}
\end{theorem}


\begin{thebibliography}{999}
\bibitem{chatzidakis} Z. Chatzidakis, `Model theory of finite fields and pseudo-finite fields', Ann. Pure Appl. Logic  88  (1997),  95--108.
\bibitem{chatznd} Z. Chatzidakis, `Model theory of difference fields', {\em The Notre Dame Lectures}. Lecture Notes in Logic 18, Assoc. Symb. Logic Urbana, IL, 2005, pp. 45-96. 
\bibitem{cdm} Z. Chatzidakis, L. van den Dries, A.J. Macintyre, `Definable sets over finite fields', J. Reine Angew. Math. 427 (1992), 107--135.
\bibitem{ch} G. Cherlin, E. Hrushovski, {\em Finite structure with few types}, Ann. Math. Studies, Princeton University Press, Princeton, 2003.
\bibitem{elwes} R. Elwes, `Asymptotic classes of finite structures', J. Symb. Logic  72 (2007), 418--438.
\bibitem{em} R. Elwes, H.D. Macpherson, `A survey of asymptotic classes and measurable structures', in
 {\em Model Theory with Applications to Algebra and Analysis, Vol. 2} (Eds. Z. Chatzidakis, H.D. Macpherson, 
A. Pillay, A.J. Wilkie), London Math. Soc. Lecture Notes 350, Cambridge University Press, 2008, pp. 125--159.
\bibitem{gowers} T. Gowers, `Quasirandom groups', Combinatorics, Probability and Computing 17 no. 3 (2008), 363--387.
\bibitem{granger} N. Granger, `Stability, simplicity and the model theory of bilinear forms', PhD thesis, University of Manchester, 1999.
www.maths.manchester.ac.uk/~mprest/
\bibitem{hodges} W. Hodges, {\em Model Theory}, Cambridge University Press, Cambridge, 1993.
\bibitem{hrush} E. Hrushovski, `Unimodular minimal structures', J. Symb. Logic 46 (1992), 385--396.
\bibitem{h1} E. Hrushovski, `Stable group theory and approximate subgroups', J. Amer. Math. Soc. 25 (2012), 189--243.
\bibitem{Hru} E. Hrushovski. \emph{On Pseudo-Finite Dimensions.} Notre Dame Journal of Formal Logic. Volume 54 (2013), no. 3-4, 463--495.
\bibitem{hrushfrob} E. Hrushovski, `The elementary theory of the Frobenius automorphisms', arxiv:math/0406514.
\bibitem{hp0} E. Hrushovski, A. Pillay, `Weakly normal groups', in {\em Logic Colloquium  85} (ed. Paris Logic Group), 1987, North-Holland, Amsterdam,  233--244.
\bibitem{hp} E. Hrushovski, A. Pillay, `On NIP and invariant measures', J. Euro. Math. Soc. 13 (2011), 1003--1061.
\bibitem{hw} E. Hrushovski, F. Wagner, `Counting and dimensions', {\em Model Theory with Applications to Algebra and Analysis, Vol. 2} (Eds. Z. Chatzidakis, H.D. Macpherson, A. Pillay, A.J. Wilkie, Cambridge University Press, Cambridge, 2008, 161--176.
\bibitem{kestner} C. Kestner, A. Pillay, `Remarks on unimodularity', J. Symb. Logic 76 (4) (2011), 1453--1458. 
\bibitem{kiefe} C. Kiefe, `Sets definable over finite fields: their zeta functions', Trans. Amer. Math. Soc. 223 (1976), 45-59.

\bibitem{kp} B. Kim, A. Pillay, `From stability to simplicity', Bull. Symb. Logic 4 (1998) no. 1, 17-36.
\bibitem{kp2}B. Kim, A. Pillay, `Simple theories', Ann. Pure Appl. Logic 88 (1997), 149--164.
\bibitem{kuzichev} A.A. Kuzichev, `Elimination of quantifiers over vectors in some theories of vector spaces', Zeit. f\"ur Math. Logik und Grundlag. der Math., 38 (1992), 575--577.
\bibitem{landau} E. Landau, `\"Uber die Klassenzahl der bin\"aren quadratischen Formen von negativer Diskriminante', Math. Ann. 56 (1903), 671--676.
\bibitem{lp} M. Larsen, R. Pink, `Finite subgroups of algebraic groups', J. Amer. Math. Soc. 24 (2011), 1105--1158.
\bibitem{lst} M. Larsen, A. Shalev, P. Tiep, `The Waring problem for finite simple groups', Ann. Math.  174 (2011), 1885--1950.
\bibitem{ms} H.D. Macpherson, C. Steinhorn, `One-dimensional Asymptotic Classes and
Measurable Structures', Trans. Amer. Math. Soc. 360 (2008), 411--448.
\bibitem{mt0} H.D. Macpherson, K. Tent, `Stable pseudofinite groups', J. Alg. 312 (2007), 550--561.
\bibitem{mt} H.D. Macpherson, K. Tent, `Pseudofinite groups with NIP theory and definability in finite simple groups', in {\em Groups and model theory}, Contemp. Math. 576, Amer. Math Soc., Providence, 2012, pp. 255-267. 
\bibitem{neumann} B.H. Neumann, `Groups covered by permutable subsets', J. London Math. Soc. 29 (1954), 236--248. 
\bibitem{np} N. Nikolov, L. Pyber, `Product decompositions of quasirandom groups and a Jordan type theorem', J. Euro. Math. Soc. 13 (2011), 1063--1077
\bibitem{pierce} D. Pierce, `Model-theory for vector-spaces over unspecified fields', Arch. Math. Logic, 48 (2009), 421--436.
\bibitem{pillay98} A. Pillay, `Definability and definable groups in simple theories',  J. Symb. Logic  63 (1998), 788--796. 
\bibitem{pillay1} A. Pillay, P. Tanovi\'c, `Generic stability, regularity, and quasiminimality', in {\em Models, logics, and higher-dimensional categories}, CRM Proc. Lecture Notes 53, Amer. Math. Soc., Providence, 2011, pp. 189--211.
\bibitem{psw} A. Pillay, T. Scanlon, F. Wagner, `Supersimple fields and division rings', Math. Research Letters 5 (1998), 473-483. 
\bibitem{ps} A. Pillay, S. Starchenko, `Remarks on Tao's algebraic regularity lemma', arXiv1310.7538.
\bibitem{ryten} M.J. Ryten, {\em Model theory of finite difference fields and simple groups}, PhD thesis, University of Leeds, 2007, http://www1.maths.leeds.ac.uk/Pure/staff/macpherson/ryten1.pdf
\bibitem{tao} T. Tao, `Expanding polynomials over finite fields of large characteristic, and a regularity lemma for definable sets', arXix:1211.2894
\bibitem{wagner} F.O. Wagner, {\em Simple theories}, Kluwer, Dordrecht, 2000.
\end{thebibliography}
\end{document}